\documentclass[12pt]{amsart}
\usepackage{amssymb,color,epsfig}
\usepackage{fullpage}

\newtheorem*{mainthm}{Main Theorem}
\newtheorem*{qal}{Quasi-Additivity Law}

\newtheorem*{covlem}{Covering Lemma}

\newtheorem*{bbrmt}{Beau Bounds (Refined Main Theorem)}

\newcommand{\from}{\mathpunct:}

\newtheorem{thm}{Theorem}[section]
\newtheorem{cor}[thm]{Corollary}
\newtheorem{lem}[thm]{Lemma}
\newtheorem{prop}[thm]{Proposition}
\theoremstyle{remark}
\newtheorem{rem}{Remark}[section]

\theoremstyle{definition}

\numberwithin{equation}{section}
\numberwithin{figure}{section}

\def\note#1{}

\def\sss{\subsubsection}

\newcommand{\di}{\partial}

\newcommand{\ra}{\rightarrow}

\newcommand{\imply}{\Rightarrow}

\def\ssk{\smallskip}
\def\msk{\medskip}

\def\nin{\noindent}

\def\sm{\smallsetminus}

\newcommand{\diam}{\operatorname{diam}}

\newcommand{\ed}{\operatorname{d}}

\newcommand{\cl}{\operatorname{cl}}
\newcommand{\inter}{\operatorname{int}}
\renewcommand{\mod}{\operatorname{mod}}

\newcommand{\orb}{\operatorname{orb}}

\newcommand{\id}{\operatorname{id}}

\newcommand{\depth}{\operatorname{depth}}
\newcommand{\bidepth}{\operatorname{bidepth}}

\renewcommand{\d}{{\diamond}}

\newcommand{\eps}{{\varepsilon}}

\newcommand{\De}{{\Delta}}
\newcommand{\de}{{\delta}}
\newcommand{\la}{{\lambda}}
\newcommand{\La}{{\Lambda}}
\newcommand{\si}{{\sigma}}

\newcommand{\Om}{{\Omega}}
\newcommand{\om}{{\omega}}

\newcommand{\al}{{\alpha}}

\newcommand{\bk}{{\boldsymbol{\kappa}}}

\newcommand{\AAA}{{\mathcal A}}
\newcommand{\BB}{{\mathcal B}}
\newcommand{\CC}{{\mathcal C}}

\newcommand{\EE}{{\mathcal E}}

\newcommand{\GG}{{\mathcal G}}

\newcommand{\HH}{{\mathcal H}}
\newcommand{\KK}{{\mathcal K}}
\newcommand{\LL}{{\mathcal L}}
\newcommand{\MM}{{\mathcal M}}

\newcommand{\OO}{{\mathcal O}}

\newcommand{\RR}{{\mathcal R}}
\newcommand{\TT}{{\mathcal T}}

\newcommand{\WW}{{\mathcal W}}

\newcommand{\YY}{{\mathcal Y}}

\newcommand{\A}{{\Bbb A}}
\newcommand{\C}{{\Bbb C}}

\newcommand{\D}{{\Bbb D}}

\newcommand{\N}{{\Bbb N}}
\newcommand{\Q}{{\Bbb Q}}
\newcommand{\R}{{\Bbb R}}
\newcommand{\T}{{\Bbb T}}

\newcommand{\Z}{{\Bbb Z}}

\newcommand{\f}{{\bf f}}
\newcommand{\g}{{\bf g}}
\newcommand{\h}{{\bf h}}
\renewcommand{\i}{{\mathbf{ i}}}

\newcommand{\n}{{\mathbf{n}}}
\newcommand{\q}{{\mathbf{q}}}
\newcommand{\p}{{\mathbf{p}}}
\newcommand{\psimod}{{ \mathbf{mod}\, }}

\def\Bf{{\mathbf{f}}}

\def\Bd{{\mathbf{d}}}

\def\BE{{\mathbf{E}}}

\def\Bpsi{{\boldsymbol{\Psi}}}

\def\B0{{\mathbf{0}}}
\def\BUps{{\boldsymbol{\Upsilon}}}
\def\BF{{\mathbf{F}}}
\def\BLa{{\boldsymbol{\La}}}
\def\BP{{\mathbf{P}}}
\def\BQ{{\mathbf{Q}}}

\def\BV{{\mathbf{V}}}
\def\BU{{\mathbf{U}}}
\newcommand{\BW}{{\mathbf{W}}}
\def\BY{{\mathbf{Y}}}
\def\BZ{{\mathbf{Z}}}
\def\BE{{\mathbf{E}}}

\def\Ups{{\Upsilon}}
\def\BUps{{\boldsymbol{\Upsilon}}}
\def\BLa{{\boldsymbol{\La}}}

\newcommand{\area}{\operatorname{area}}

\newcommand{\fjset}{FJ-set}

\catcode`\@=12

\def\Empty{}
\newcommand\oplabel[1]{
  \def\OpArg{#1} \ifx \OpArg\Empty {} \else
  	\label{#1}
  \fi}
		
\newcommand{\comm}[1]{}
\newcommand{\comment}[1]{}

\def\makeabbrevs{%
\def\o{\omega}\def\g{\gamma}\def\G{\Gamma}\def\h{\hat}\def\d{\delta}\def\D{\Delta}%
\def\O{\Omega}\def\b{\beta}\def\l{\lambda}}

\begin{document}

\bigskip\bigskip

\title[Decorations]{A priori bounds for some\\ infinitely renormalizable quadratics:
  {\small II. Decorations.} }
\author {Jeremy Kahn and Mikhail Lyubich}

\begin{abstract}
A decoration of the Mandelbrot set $M$ is a part of $M$ cut off by two external rays
landing at some tip of a satellite copy of $M$ attached to the main cardioid.
In this paper we consider infinitely renormalizable quadratic polynomials satisfying the decoration condition,
which means that the combinatorics of the renormalization operators involved is selected from
a finite family of decorations. For this class of maps we prove {\it a priori} bounds.
They imply local connectivity of the corresponding Julia sets and the Mandelbrot set
at the corresponding parameter values.
\end{abstract}

\setcounter{tocdepth}{1}

\maketitle

\thispagestyle{empty}
\def\IMSmarkvadjust{0 pt}
\def\IMSmarkhadjust{0 pt}
\def\IMSmarkhpadding{0 pt}
\def\IMSpubltext{Published in modified form:}
\def\SBIMSMark#1#2#3{
 \font\SBF=cmss10 at 10 true pt
 \font\SBI=cmssi10 at 10 true pt
 \setbox0=\hbox{\SBF \hbox to \IMSmarkhpadding{\relax}
                Stony Brook IMS Preprint \##1}
 \setbox2=\hbox to \wd0{\hfil \SBI #2}
 \setbox4=\hbox to \wd0{\hfil \SBI #3}
 \setbox6=\hbox to \wd0{\hss
             \vbox{\hsize=\wd0 \parskip=0pt \baselineskip=10 true pt
                   \copy0 \break%
                   \copy2 \break%
                   \copy4 \break}}
 \dimen0=\ht6   \advance\dimen0 by \vsize \advance\dimen0 by 8 true pt
                \advance\dimen0 by -\pagetotal
	        \advance\dimen0 by \IMSmarkvadjust
 \dimen2=\hsize \advance\dimen2 by .25 true in
	        \advance\dimen2 by \IMSmarkhadjust

%
%
  \openin2=publishd.tex
  \ifeof2\setbox0=\hbox to 0pt{}
  \else 
     \setbox0=\hbox to 3.1 true in{
                \vbox to \ht6{\hsize=3 true in \parskip=0pt  \noindent  
                {\SBI \IMSpubltext}\hfil\break
                \input publishd.tex 
                \vfill}}
  \fi
  \closein2
  \ht0=0pt \dp0=0pt
 \ht6=0pt \dp6=0pt
 \setbox8=\vbox to \dimen0{\vfill \hbox to \dimen2{\copy0 \hss \copy6}}
 \ht8=0pt \dp8=0pt \wd8=0pt
 \copy8
 \message{*** Stony Brook IMS Preprint #1, #2. #3 ***}
}

\SBIMSMark{2006/6}{August 2006}{}

\tableofcontents

\section{Introduction}

\begin{figure}
\centerline{\includegraphics[width=12cm]{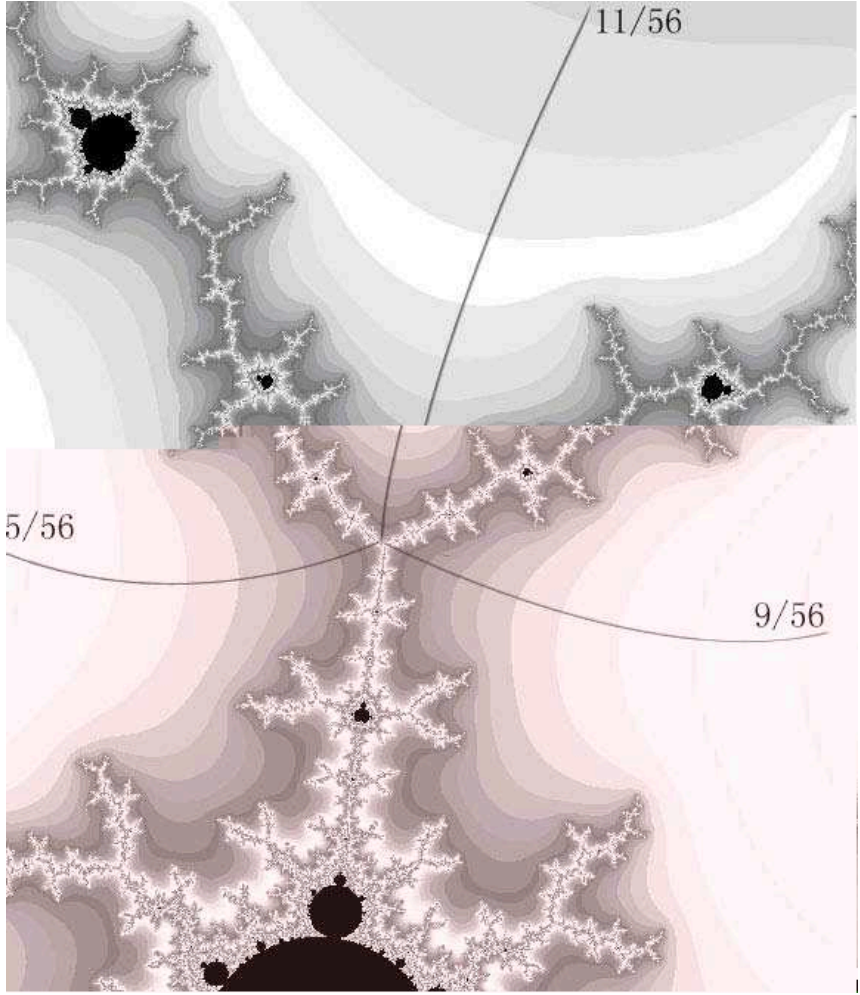}}
\caption{Two decorations of the Mandelbrot set.
   The little Mandelbrot sets inside specify renormalization combinatorics of type $(3,1)$.  \label{limbs}}
 \end{figure}

  A {\it decoration} of the Mandelbrot set $M$ (called also a {\it Misiurewicz limb})
$\LL$ is a part of $M$ cut off by two external rays
 landing at some tip of a satellite copy of $M$ attached to the main cardioid,
see Figure \ref{limbs} (see \S \ref{sec: first king} for the precise dynamical definition).
In this paper we consider infinitely
renormalizable quadratic polynomials satisfying the decoration condition,
which means that the combinatorics of the renormalization operators involved is selected from
a finite family of decorations $\LL_k$. (For instance, real infinitely renormalizable maps
satisfy a decoration condition if and only of non of the renormalizations is of doubling type.)

An infinitely renormalizable quadratic map $f$ is said to have {\it a priori bounds}
if its renormalizations can be  represented by quadratic-like maps  $R^n f: U_n\ra V_n$
with $\mod (V_n\sm U_n)\geq \eps >0$, $n=1,2\dots$.

Our goal is to prove the following result:

\begin{mainthm}
  Infinitely renormalizable quadratic maps satisfying the decoration condition have
a priori bounds.
\end{mainthm}

By \cite{puzzle}, this implies:

\begin{cor}\label{JLC and MLC}
 Let $f_c: z\mapsto z^2+c$ be an infinitely renormalizable quadratic map satisfying the decoration condition.
Then the Julia set $J(f_c)$ is locally connected, and the Mandelbrot set $M$ is locally connected at $c$.
\end{cor}

In this paper we will deal only with the case of sufficiently high periods:

\begin{thm}\label{high periods thm}
 Given finitely many decorations $\LL_k$,
 there exists a $\underline p$ such that
any infinitely renormalizable quadratic map satisfying the decoration condition
with decorations $\LL_k$ and renormalization periods $p\geq \underline p$ has
a priori bounds.
\end{thm}

The complementary case of ``bounded combinatorics'' is dealt in \cite{K}.

\begin{rem}
Theorem \ref{high periods thm} sounds similar in spirit to the {\it a priori} bounds
  of \cite{puzzle}. However, the   ``high type'' condition of \cite{puzzle} is stronger then the
  above high period condition, while the ``secondary limb condition'' of \cite{puzzle} is weaker than
the decoration condition. Also, our proof of Theorem~\ref{high periods thm} is compatible with the
  proof of \cite{K}, so that they can be combined into the Main Theorem.
\end{rem}

Let us now outline the structure of the paper.

In the next section, \S \ref{principal nest}, we will describe a necessary combinatorial
set-up in the framework of the Yoccoz puzzle. Besides a well-known material,
it includes the construction of the {\it modified principal nest} from \cite{KL-high}
needed for dealing with maps of ``high type''.


In \S \ref{pseudo-puzzle} we summarize necessary information about
{\it pseudo-quadratic-like maps} defined in \cite{K},
and introduce a {\it pseudo-puzzle} by applying the ``pseudo-functor'' to the puzzle.
In this way we make domains of the return maps more canonical,
which spares us from the need to control geometry of  external rays.

From now on, the usual puzzle will serve only as a combinatorial frame,
while all the geometric estimates will be made on the pseudo-puzzle.
This is needed for this paper per se, as well as for making connection to the
case of bounded combinatorics \cite{K}.
 Only at the last moment (\S \ref{conclusion})
we return back to the standard quadratic-like context.

In \S \ref{covering lemma sec} we formulate the analytical results of \cite{covering lemma},
the Quasi-Additivity Law and the Covering Lemma, in the pseudo context.
They will be our main analytical tools.

In \S \ref{moduli} we prove the main results of the paper.
To prove {\it a priori} bounds,  we show  that
if some renormalization has a small modulus, then this modulus will improve on some deeper level.
The main place where the decoration condition plays the role is on the top of the puzzle,
when we compare the modulus of the first annulus of the pseudo-puzzle to the modulus of the original
pseudo-quadratic-like map.

\begin{rem}
   Strictly speaking, bounded combinatorics treated in \cite{K}
and high combinatorics treated by Theorem \ref{high periods thm} do not cover
the oscillating combinatorial types. However, these theorems follow from results
on moduli improving that together cover everything.
\end{rem}


\begin{rem}
Our proof of {\it a priori bounds} (Main Theorem) applies without changes in the case
of unicritical maps of higher degree. However, the proof of MLC at the corresponding parameters
(Corollary~\ref{JLC and MLC}) given in \cite{puzzle} exploits some special geometric features of
quadratic maps. In \cite{Ch} part of \cite{puzzle} is combined with a new method developed in
\cite{AKLS} to prove Corollary \ref{JLC and MLC} in the higher degree case as well.
\end{rem}

\subsection{Terminology and Notation}
  $\N =\{1,2,\dots\}$ is the set of natural numbers; $\Z_{\geq 0} = \N \cup \{0\}$;
$\D=\{z:\,  |z|<1\}$ is the unit disk, and $\T$ is the unit circle;\\
$\A(r,R)= \{z: r<|z|<R\}$ is the annulus of modulus $\frac{1}{2\pi} \log (R/r)$;\\
$\Pi(h)= \{z| 0< \Im z<h\}$ is the horizontal strip.

  A  {\it topological disk} means a simply connected  domain in some Riemann surface $S$.
A {\it continuum} $K$ is a connected closed subset in $S$. 
It is called {\it full} if all components of $S\sm K$ are unbounded.
We say a subset $K$ of a plane is an \fjset\ (for ``filled Julia set'')
if $K$ is compact, connected, and full.

  We let $\orb(z)\equiv \orb_g(z)= (g^n z)_{n=0}^\infty$ be the {\it orbit} of $z$ under a map $g$.

  Given a map $g: U\ra V$ and an open topological disk $D\subset V$,
components of $g^{-1}(D)$ are called {\it pullbacks}  of $D$ under $g$.
If the disk $D$ is closed, we define pullbacks of $D$ as the closures of the pullbacks of $\inter D$.%
\footnote{Note that the pullbacks of a closed disk $D$ can touch one another,
so they are not necessarily connected components of $g^{-1}(D)$.}
In either case,
given a connected set $X\subset g^{-1}(\inter D)$, we let $g^{-1} (D)|X$ be the pullback of $D$
containing $X$.

We let $ x\oplus y = (x^{-1} + y^{-1})^{-1}$ be the {\it harmonic sum} of $x$ and $y$
(it is conjugate to the ordinary sum by the inversion map $x\mapsto x^{-1}$).
Similarly,  $ x\ominus y = (x^{-1} - y^{-1})^{-1}$ stands for the  harmonic difference.

%

\subsection{Acknowledgement}
  We thank Tao Li for making Figure \ref{limbs}.
This work has been partially supported by the NSF, NSERC, the Guggenheim and  Simons Foundations.
Part of it was done during the authors' visit to the IMS at Stony Brook and the Fields Institute in Toronto.
We are thankful to all these Institutions and Foundations.

\section{Yoccoz puzzle, decorations, and the Modified Principal Nest}\label{principal nest}

Let $(f_\la: U_\la'\ra U_\la)$ be a quadratic-like family over a disk $\La\subset \C$.
Assume that this family is good enough (proper and unfolded), so that the associated
Mandelbrot set $M=M(f_\la)$ is canonically homeomorphic to the standard Mandelbrot set (see \cite{DH}).
In fact, most of the time we will be dealing with a single map $f=f_\la$ from our family,
so that we will usually suppress the label $\la$  in the notation.
(We need a one parameter family only to introduce different combinatorial types of the maps
under consideration.)

We assume that the domains $U'$ and $U$ are smooth disks,
$f$ is even,
and we  normalize $f$ so that $0$ is its critical point.

We let $U^m = f^{-m}(U)$. The boundary of $U^m$ is called
the {\it equipotential of level $m$}.

\subsection{Top of the Yoccoz puzzle and  decorations}\label{sec: first king}
 By means of straightening, we can define external rays for $f$. They form a foliation of $V\sm K(f)$
orthogonal  to the equipotential $\di U$.
 The map $f$ has one {\it non-dividing} fixed point $\beta$ (landing point of the external ray
with angle $0$), and one {\it dividing} fixed point $\alpha$. There are $\q>1$ external rays $\RR_i$ landing at
$\alpha$ which are cyclically permuted by the dynamics with rotation number $\p/\q$, see \cite{M-rays}
($\p/\q$ is also called the {\it combinatorial rotation number} of $\alpha$).
These rays divide $U$ into $\q$ (closed) topological disks $Y^0_i$ called the {\it Yoccoz puzzle pieces} of depth 0.
Let $Y^0\equiv Y^0_0$ stand for the critical puzzle piece, i.e., the one containing 0.

Let us consider  $2\q$  rays of  $f^{-1} (\cup\RR_i)$.
They divide $U'$ into $2\q-1$ (closed) disks called Yoccoz puzzle pieces of depth 1.
Let $Y^1$ stand for the critical puzzle piece of depth 1.
There are also $\q-1$  puzzle pieces $Y^1_i$ of depth 1 contained in the corresponding off-critical pieces of depth 0.
All other puzzle pieces of depth 1 will be denoted $Z_i^1$.
They are attached to the symmetric point $\alpha'=-\alpha$.

The puzzle pieces will be labeled in such a way that $f(Y_i^1) =  Y^0_{i+1}$,
$i=0,\dots, \q-1$, and $Z^1_i= - Y^1_i$.
We let
$$
   L= \bigcup_{i=1}^{\q-1} Y_i^1; \quad  R= -L= \bigcup_{i=1}^{\q-1} Z_i^1.
$$

Puzzle pieces $Y^m_j$ of depth $m$ are pullbacks of $f^{-m}(Y^0_i)$.
They tile the neighborhood of $K(f)$ bounded by the equipotential $\di U^m $.
Each of them is bounded by finitely many arcs of this equipotential and finitely many
external rays of $f^{-m}(\RR_i)$.
If $f^m(0)\not=\alpha$, then  there is one puzzle piece of depth $m$ that contains the critical point $0$.
It is called {\it critical} and is labeled as $Y^m\equiv Y^m_0$.
These pieces are nested around the origin:
$$
   Y^0\supset Y^1\supset Y^2 \dots \ni 0.
$$

Let us consider a puzzle piece $Y=Y^m_i$.
Different arcs of $\di Y$  meet at the {\it corners} of $Y$.
The corners where two external rays meet will be called {\it vertices} of $Y$;
they are $f^m$-preimages of $\alpha$.
Let $K_Y= K(f)\cap Y$. It is a closed connected set that meets the boundary $\di Y$
at its  vertices. Moreover, the external rays meeting at a vertex $v\in \di Y$ chop off from
$K(f)$ a continuum $S_Y^v$, the component of $K(f)\sm \inter Y$ containing $v$.

The critical value $f^q(0)$ belongs to the puzzle piece $Y^0$.
If in fact it belongs to $Y^1$ then the map $Y^{\q+1}\ra Y^1$ is a double branched covering.
It is not a quadratic-like map, though, since the boundaries of $Y^1$ and $Y^{\q+1}$ overlap over four external rays
landing at $\alpha$ and $\alpha'$. However,  by slight ``thickening'' of the domain of this map (see \cite{M}),
it can be turned into a quadratic-like map $g$ such that
$$
     K(g) = \{ z:\, f^{\q m}z\in Y^1, \, m=0,1,2,\dots \}.
$$

  The map $f$ is called   {\it satellite renormalizable}
(or, {\it immediately renormalizable})  if  the Julia set $K(g)$ is connected, i.e.,
if the critical point never escapes $Y^1$:
 $$
           f^{\q m}(0)\in Y^1,\quad  m = 0,1,2\dots.
$$
The set of immediately renormalizable parameter values (with a given combinatorial rotation number $\p/\q$)
assemble a {\it satellite copy $M_{\p/\q}$ of $M$} attached to the main cardioid at the parabolic point with rotation
number $\p/\q$.  The parameters $t\in M_{\p/\q}$ for which the critical point eventually lands at $\alpha$
(i.e., $f_t^{\q n} = \alpha'$ for some $n\in \N$) are called the {\it tips} of $M_{\p/\q}$.

\ssk
If $f$ is not satellite renormalizable, then
there exists an $n\in \N$ such that $f^{\q\n}(0)$ belongs
to some puzzle piece $\inter Z_\bk^1$.
Let $\n$ be the smallest such $n$.
In this case, we let
$$
   V^0= f^{-\n\q}(Z_\bk^1)|0= Y^{\n\q+1}.
$$

Each puzzle piece $Z^1_j$ has $2^m$ univalent pullbacks under the $2^m$-covering $f^{\q m}: Y^{\q m}\ra Y^0$,
$m=1,\dots, \n-1$. We label these pullbacks (for all $j$) as $Z^{1+\q m}_i$.  Then
\begin{equation}\label{escape route}
    f^{\q m}(0) \in Z^{\q(\n-m)+1}_{\bk_m}, \quad m=1,\dots, \n,
\end{equation}
for some sequence $\bar\kappa = (\kappa_1,\dots, \kappa_n=\bk)$ called the
{\it escape route} of the critical point. The escape route specifies
the tip $t=t_{\kappa_1\dots\kappa_{\n-1}}$ of $M$ such that $f_t$ satisfies
(\ref{escape route}) for $m<\n$, while $f^{\q\n}=\alpha'$.

There are $\q$ parameter rays landing at each tip $t$ of $M_{\p/\q}$.
They chop off $\q-1$ {\it decorations} $\LL_{\bar\kappa}$
(the components of $M\sm \{t\}$ that do not intersect the main cardioid) from $M$.
The limb $\LL_{\bar\kappa}$ attached to $t$ is specified by the puzzle piece $Z^1_\bk$ containing $f^{\q\n}(0)$.
Note that there are only finitely many decorations with bounded $\q$ and $\n$.

Let $P=Y^{(\n-1)\q+1}$. The piece $P$ has $2^\n$ vertices each of which is a preimage of $\alpha$
of some depth $\q\ m$ with $m\leq \n$ (and it takes into account {\it all} preimages of $\alpha$ in $Y^1$
 up to depth  $\q\n$).
\begin{figure}[htbp]
\begin{center}
\makeabbrevs
\input{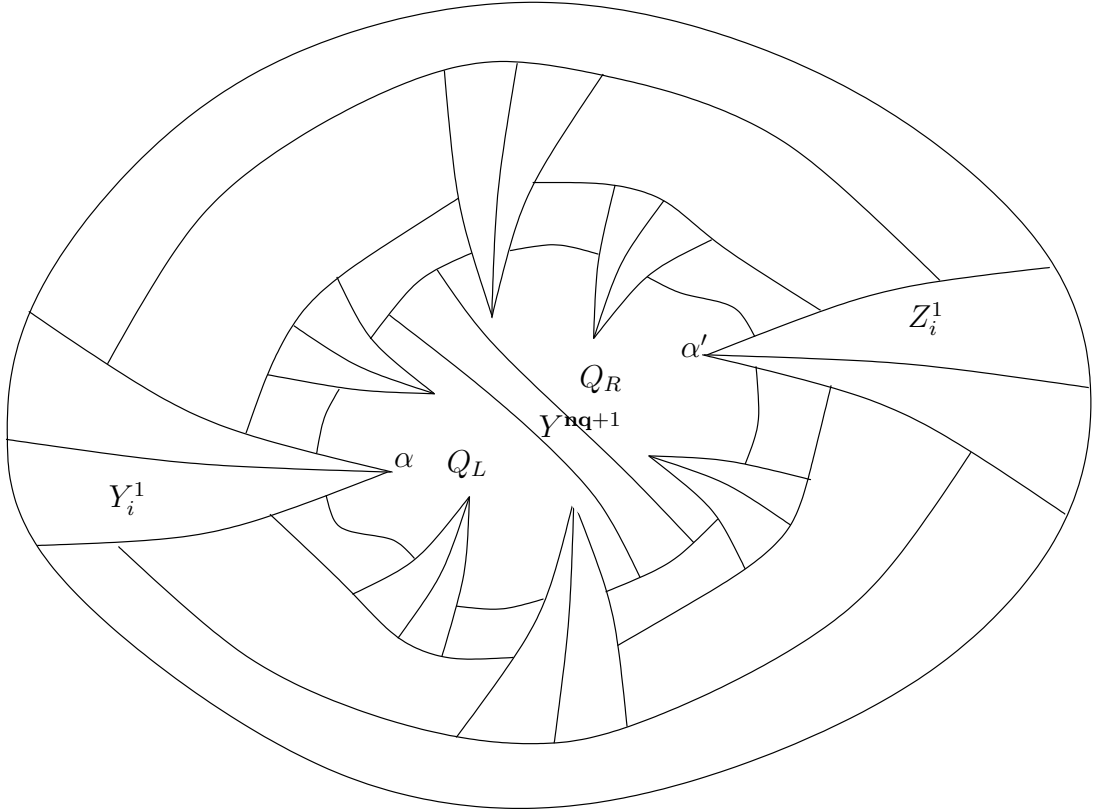}
\caption{The top of the Yoccoz puzzle}
\label{puzzle}
\end{center}
\end{figure}

\ssk
{\it Remark 2.1.}  Any $f^{\q m}$-preimage  $t\in Y^1$ of $\alpha$, $m\leq \n$,
   can be naturally labeled by a dyadic number
  $i/ 2^m\in \Q/\Z$ (with odd $i$). Here $\alpha$ is labeled by $0$, $\alpha'$ is labeled by $1/2$,
  and in general, the dyadic expansion of $i$ is  $0.\,\eps_1\dots\eps_m\, 1$,
  where $\eps_l$ is equal $+1$ or $-1$ depending on whether the rays landing at $f^{\q(l-1))}(t)$
  are ``below'' the chord connecting $\alpha$ and $\alpha'$, or above it.%
 Then the pullbacks of $Z^j$ attached to $t$ can be labeled as $Z_j(i/ 2^m)$, $j=1,\dots, \q-1$.
  A decoration assumes labeling $Z_j(i / 2^\n) $ if $f^\q(0)\in Z_j(i/2^\n)$.

\medskip
Note that $f^\q(P)\supset P$ and  the critical value  $f^\q(0)$
does not belong to $P$. Hence $P$ has two univalent $f^\q$-pullbacks,
$Q_L$ and $Q_R$ (of depth $\q \n+1$),  inside $P$.
The puzzle piece $Q_L$ is attached to the fixed point $\alpha$
 while $Q_R$ is attached to $\alpha'$. Each of them shares two external rays with $V^0$.

\begin{lem}\label{separating family}
For any  vertex $v$ of $P$,
there exists a puzzle piece $Q^v\subset P$ of depth $(2\n-1)\q+1$
attached to the boundary rays of $P$ landing at $v$,
 which is a univalent $f^{\n\q}$-pullback of $P$.
Moreover, these puzzle pieces are pairwise disjoint.
\end{lem}

\begin{proof}
Let $g= f^\q |\, Q_L\cup Q_R$.
The domain of $g^\n$ consists of $2^\n$ components each of which is a univalent pullback of $P$.
Each of these components contains a single $g^\n$-preimage $v$ of $\alpha$,
and is attached to the pair of the boundary external rays of $P$ landing at $v$.
This is the desired puzzle piece $Q^v$.
\end{proof}

\comm{*****
Let $a,b\in K$.
We say that a puzzle piece $Q$ separates
$a$ from $b$
if $a$ and $b$ belong to different components of $K\sm K_Q$.
We say that a family of puzzle pieces $Q_i$ {\it doubly separates} a family of
points $a_j\in K$, if any two points, $a_j$ and $a_k$, are separated
by at least two puzzle piece $Q_i$.

\begin{lem}\label{separating family}
   There is a family $Q_i\subset Q_L\cup Q_R$
of pairwise disjoint univalent pullbacks of
the puzzle piece $P$ under the map $f^{\q}: Q_L\cup Q_R\ra P$
that doubly separates all the vertices of $P$. Moreover, $\depth Q_i\leq d(\q,\n).$
\end{lem}
***************}

Given two vertices, $v$ and $w$, of $P$,
we let $T^{vw}_P= K_P\sm (Q^v\cup Q^w)$.
Notice that $T^{vw}$ {\it separates} $v$ from $w$ in the sense that $v$ and $w$ belong to different components
of $K_P\sm T_{vw}$.

\subsection{Modified principle nest}
  Given a critical puzzle piece $V$, 
let us consider the first return $f^l 0$, $l\in \N$,  of the critical point to $V$
(whenever it exists).
The corresponding pullback $W= f^{-l}(V)|0$ of $V$ is called
the {\it central domain of the first return map to $V$},
or briefly, the {\it first child} of $V$.
Under these circumstances, $W\subset V$ and the first return map $f^l: W\ra V$ is a double branched covering.

Under the above circumstances,
we also consider the the first moment $k\in \N$ such that $f^{kl}(0)\not\in W$
and then the  first return $f^t(0)\in W$, $t>kl$, back to $W$ (whenever these moments are well defined).
We call it the {\it fine return} to $W$, and the corresponding pullback $A=f^{-t}(W)|0$ the {\it fine child} of $W$.
The map $f^t: A\ra W$ is a double branched covering. Note that if  $f^l(0)\not\in W$,
the fine return coincides with the first return.

In \cite{KL-high} we have constructed a {\it (Modified) Principle Nest} of critical puzzle pieces
$$
  E^0\Supset E^1\Supset E^1\Supset\dots\Supset E^{\chi-1}\Supset E^\chi
$$
and corresponding quadratic-like maps $g_n: E^n\ra E^{n-1}$. 
Here for odd $n$,  $E^n$ is the first child of $E^{n-1}$ and $g_n: E^n\ra E^{n-1}$ is the corresponding  first return map.
For odd $n$, $E^n$ is the fine child of $E^{n-1}$ and $g_n: E^n\ra E^{n-1}$ is the corresponding fine return map.
We let $g\equiv g_0$.

If the map $f$ is renormalizable then  the Principle Nest terminates at some odd level $\chi$.
In this case, the last quadratic-like map $g_\chi: E^\chi\ra E^{\chi-1}$  has connected Julia set
and represents the {\it primitive renormalization} $Rf$ of $f$.
The renormalization level $\chi$ is also called the {\it height} of the nest.

Primitively renormalizable parameter values assemble a maximal {\it primitive copy} $M'$ of
the Mandelbrot set $M$. This copy specifies
the {\it combinatorics} of the renormalization in question.
In particular, it determines the parameters $\q$, $\n$, the height $\chi$,
and the renormalization period $p$.

{\it In what follows we will assume that $f$ is primitively renormalizable.}
We let $\KK=K(\RR f)$ be the {\it little (filled) Julia set } of $f$,
and we let $p$ be the renormalization period, i.e., $g_\chi= f^p$ so that $f^p(\KK)=\KK$.
We let $\KK_i= f^i(\KK)$, where $i$ is taken mod $p$, which are also called ``little Julia sets''.


\msk
It is important to note that the maps $g_n$ admit analytic {\it extensions} $\tilde E^n\ra \hat E^{n-1}$
such that $E^n\subset \tilde E^n \subset \hat E^n\subset E^{n-1}$ and for odd $n$,
$\hat E^n=E^{n-1}$ \cite{KL-high}, \S 2.4. For $n=0$, we let $\hat E^0 =\tilde E^0 = Y^{\q\n}$.
Then $ f^{\q\n}: \tilde E^0 \ra Y^0$ is a branched covering of degree $2^\n$.

\msk
The following useful observation will be used many times:

\begin{lem}[Telescope]\label{telescope}
  Let $X_k$ be a sequence of topological disks, $k=0,1,\dots, m$,
and let $\phi_k : X_k \ra \phi(X^k)$  be  branched coverings of degree $d_k$ such that
$\phi(X_k)\supset X_{k+1}$. Let $\Phi= \phi_{n-1}\circ\dots\circ \phi_0$ (wherever it is defined),
and let $P\subset X_0$ be a component of its domain of definition.
Then $\Phi: P \ra V_n$ is a branched covering of degree at most $d_0\cdot\cdot\cdot d_{n-1}$.
\end{lem}

If the renormalization $Rf$ is also renormalizable then $f$ is called twice renormalizable,
and $R^2 f$ stands for its second renormalization. Proceeding this way, we can define
{\it infinitely renormalizable} maps $f$, and let $R^n f$ be their $n$-fold renormalizations.
The {\it combinatorics} of an infinitely renormalizable map is a sequence
of little Mandelbrot copies $M^{(n)}$ that determine the combinatorics of the renormalizations $R^n f$.
It determines the sequence of the parameters $\q_n$, $\n_n$, the heights $\chi_n$,
and the periods $p_n$ of the corresponding renormalizations.

We say that an infinitely renormalizable $f$ satisfies the {\it decoration condition}
if all the little copies $M^{(n)}$ belong to  finitely many decorations $\LL_k$.
Equivalently, {\it the parameters $\q_n$ and $\n_n$ are bounded}.



\subsection{Geometric puzzle pieces}

In what follows  we will deal with more general puzzle pieces.

Given  a puzzle piece $Y^m_i$, let $Y^{m,l}_i$ stand for a Jordan disk bounded by the same external rays
as $Y^m_j$ and arcs of equipotentials of level $l$ (so $Y^{m,m}_i=Y^m_i$).
Such a  disk will be called a puzzle piece of {\it bidepth} $(m,l)$.

A {\it geometric puzzle piece} of bidepth $(m,l)$ is a closed Jordan domain which is
the union of several puzzle pieces of the same bidepth.
As for ordinary pieces, a pullback of a geometric puzzle piece of bidepth $(m,l)$
under some iterate $f^k$ is a geometric puzzle piece of bidepth $(m+k, l+k)$.
Note also that if $P$ and $P'$ are geometric puzzle pieces with%
\footnote{the inequality between bidepths is understood componentwise}
$\bidepth P \geq \bidepth P'$
and $K_P\subset K_{P'}$ then $P\subset P'$.

The family of geometric puzzle pieces of bidepth $(m,l)$ will be called $\YY^m(l)$.
Given a geometric puzzle piece $Y\in \YY^m(l)$, we let $Y(k)$ be the puzzle piece
bounded by the same external rays as $Y$ truncated by the equipotential of level $k$.
(In particular, $Y(l)=Y$.)

\msk

Any puzzle piece $Y\in \YY^m(l)$ admits the following combinatorial representation.
Let $\theta_i$ be the cyclically ordered  angles of the external rays $\RR_i$ that bound $Y$.
Let us consider the straight rays $R_i$ in $\C\sm \D$ of angles $\theta_i$
truncated by the circle $\T_r$ of radius $r=1/2^l$.
If two consecutive rays, $\RR_i$ and $\RR_{i+1}$, land at the same vertex of $Y$,
let is connect $R_i$ to $R_{i+1}$ with a hyperbolic geodesic in $\D$.
Otherwise $\RR_i$ and $\RR_{i+1}$ are connected with an equipotential arc.
Then let us connect $R_i$ to $R_{i+1}$ with the appropriate arc of $\T_r$.
We obtain a Jordan curve that bounds the {\it combinatorial model} $M_Y$ of $Y$.

 The arcs $\om_i$ of  $\T\cap M_Y$ correspond to  the
``external arcs'' of the Julia piece $K_Y$.
 They have length $2\pi \la$, where $\la$ is called the {\it combinatorial length}
of the corresponding external arc of $K_Y$.
In case $Y$ is a dynamical puzzle piece,
all the external arcs of $Y$ have  the same combinatorial length
\begin{equation}\label{comb length}
    \frac{2^k}{(2^\q-1)2^m}, \quad k\in \{0,1,\dots, q-1\},
\end{equation}
where the choice of $k$ depends on  the puzzle piece $f^m(Y)$ of depth 0
(For instance, $k=0$ when $f^m(Y)$ contains the critical value $f(0)$,
 while $k=q-1$ when $f^m(Y)$ contains the critical point 0.)

It follows that
for a geometric puzzle piece $Y$ of depth $m$,
the combinatorial length of its external arcs is at least $2^{-(\q+m)}$.

\msk
Let us now consider a geometric puzzle piece  $Z^0=  - Y^0$ of bidepth (1,0).

\begin{lem}\label{qn-pullbacks}
  Let $z\in K(f)$,  $f^{\q\n} z \in Z^0$ and let $P = f^{-\q\n}(Z^0)|z$.
Then $P\Subset \inter Y^0$ or  $P\Subset \inter Z^0$.
\end{lem}

\begin{proof}
  $P$ is a geometric puzzle piece of bidepth $(\n\q+1, \n\q) $. But $E^0=Y^{\n\q+1}$ is a puzzle piece of
depth $\n\q+1$ such that $f^{\q\n}(E^0)=Z_\bk$, where $\inter Z_\bk \cap Z^0=\emptyset$.
It follows that $P \cap \inter E^0=\emptyset$. But $K(f)\sm \inter E^0$ consists of two $0$-symmetric
connected components $X_L\supset L\cap K(f)$ and $X_R\supset R\cap K(f)$.
We conclude that  $K_P$ is contained in one of these components, and hence it is contained in one of the
sets $K_{Z^0}$ or $K_{Y^0}$. As
 $$
    \bidepth P\geq (1,0)=\bidepth Z^0\geq (0,0)=\bidepth Y^0,
$$
 $P$ is contained in one of the
puzzle pieces $Z^0$ or $Y^0$.
\end{proof}

\subsection{Many happy returns}\label{high type}
  Here we will summarize the combinatorial construction of \cite{KL-high}, \S 1.9,
that will lead to the moduli improvement in the case of high type.

Fix an arbitrary $m$, let $N$ be the smallest even integer which is bigger than $\log_2 m + 5$,
and take any odd level $n\geq N$.
Then there exists $m/2$ returns $\La_k=g^{l_k}(E^n)$  of the domain $ E^n $ to $E^{n-N}$
with the following properties.
For any domain $\La_k$,
the map $\Psi_k= g^{l_k}: E^n \ra \La_k$ admits a holomorphic extension to a branched covering
\begin{equation}\label{3 domains}
    \Psi_k : (\Upsilon_k , \De_k, E^n) \ra (E^{n-N-1}, \La_k', \La_k)
\end{equation}
such that:

\begin{itemize}

\item [(P1)] $\deg(\Psi_k: \Upsilon_k\ra E^{n-N-1} )\leq 2^{N+m}$;

\item[(P2)] $\deg(\Psi_k: \De_k\ra \La_k')\leq d^5$;

\item [(P3)] $\Upsilon_k\subset E^{n-1}$;

\item[(P4)] There is a level $i\in [n-5, n-1]$ such that each pair of disks $(\La_k', \La_k)$ is
mapped  univalently onto  $( \hat E^i,  E^i)$ under some iterate $f^t$, $t=t(k)$;

\item[(P5)] The buffers $\La_k'\Subset E^{n-N}$ are pairwise disjoint.

\end{itemize}

\section{Pseudo-quadratic-like maps and pseudo-puzzle}\label{pseudo-puzzle}

\subsection{Pseudo-quadratic-like maps}

For a more general and detailed discussion of $\psi$-ql maps, see \cite{K}.

Suppose that $\BU'$, $\BU$ are disks,
and $i:  \BU' \to \BU$ is a holomorphic immersion,
and $f: \BU' \to \BU$ is a degree $d$ holomorphic branched cover.
Suppose further that there exist full continua $K \Subset \BU$ and $K'\Subset \BU'$
such that $K' = i^{-1}(K) = f^{-1}(K)$.
Then we say that $F=(i,f):  \BU'\ra (\BU, \BU)$  is a
$\psi$-{\it quadratic-like}
 ($\psi$-{\it ql}) map  with filled Julia set $K$.

\begin{lem}[\cite{K}]\label{ql extension}
Let $F=(i,f)\from\BU'\ra \BU$ be a $\psi$-ql map of degree $d$
 with filled Julia set $K$.
Then $i$ is an embedding in  a neighborhood of $K' \equiv f^{-1}(K)$,
and the map  $g\equiv  f\circ i^{-1}\from U'\ra U$ near $K$ is quadratic-like.

\ssk Moreover, the domains $U$ and $U'$ can be selected in such a way that
$\mod(U\sm i(U'))\geq \mu (\mod(\BU\sm K)) >0$.
\end{lem}

There is a natural $\psi$-ql map $\BU^n\ra \BU^{n-1}$,
the ``restriction'' of $(i,f)$ to $\BU^n$.
Somewhat loosely, we will use the same notation $F=(i,f)$ for this restriction.

 Let us normalize  the $\psi$-quadratic-like  maps under consideration
so that $\diam K'= \diam K=1$, both $K$ and $K'$ contain 0 and 1,
$0$ is the critical point of $f$, and $i(0)=0$.
Let us endow the space of $\psi$-quadratic-like maps (considered up to independent rescalings in the domain and the range)
with the Carath\'eodory topology.
In this topology, a sequence of normalized maps $(i_n, f_n): \BU_n'\ra\BU_n$ converges to $(i,f): \BU'\ra \BU$ if
the pointed domains $(\BU_n', 0)$ and $(\BU_n, 0)$ converge to $\BU'$ and $\BU$ respectively,
and the maps $i_n$, $f_n$ converge respectively to $i$, $f$,
uniformly on compact subsets of  $\BU'$.

\begin{lem}[compare \cite{McM}]\label{compactness}
  Let $\mu>0$. Then the space of $\psi$-PL maps $F=(i,f): \BU'\ra \BU$ such that the  Julia set $K$
is connected and  $\mod(\BU \sm K)\geq \mu$ is compact.
\end{lem}

\begin{proof}
 Let  $X_n= i^{-1}\{0,1\}$,  $Y_n= f^{-1}\{0,1\}$.
Note that both sets consist of at most 2 points and are contained in $\D$.

Then we can select a subsequence of domains $\BU_n'$, $\BU_n$ Carath\'eodory converging to some
domains $\BU'$, $\BU$, while the sets $X_n$ and $Y_n$ converge in the Hausdorff metric to some sets
$X\subset \BU'$  and $Y\subset \BU'$ that consist of at most two points and are contained in $\D$
 (we will keep the same notation for the subsequence).
Since  the maps $i_n| \BU_n'\sm X_n$ and $f_n| \BU_n'\sm Y_n$ do not assume values $0$ and $1$,
they form normal families on $\BU'\sm X $. Since these families are bounded on the sets $K_n$,
they are uniformly bounded on compact sets of $\BU'\sm (X\cup Y)$.
By the Maximum Principle, they are normal on the whole domain $\BU'$.

Let $i$ and $f$ be some limit functions of the sequences $i_n$ and $f_n$.
These functions are non-constant since they assume values $0$ and $1$.
Then $i$ is an immersion as a non-constant limit of immersions.
Also, $f: \BU'\ra \BU$ is a branched covering of degree at most 2.
Moreover,  $K'\Subset \BU$ since $\mod(\BU'\sm K')\geq \mu/2$. Hence $0\in \BU'$, and it is a critical point of $f$.
It follows that $\deg f=2$, and we are done.
\end{proof}

\subsection{Pseudo-puzzle}\label{psi-puzzle}

\subsubsection{Definitions}
   Let $(i,f) : \BU'\ra \BU$ be a $\psi$-ql map.
By Lemma~\ref{ql extension}, it admits a quadratic-like restriction $U'\ra U$ to a neighborhood of
its (filled) Julia set  $K=K_\BU$.
Here $U'$ is embedded to $U$, so we can identify $U'$ with $i(U')$
and $ f: U'\ra U $ with  $f \circ i^{-1} $.

Assume that $K$ is connected and both fixed points of $f$ are repelling.
Then we can cut $U$ by external rays landing at the $\alpha$-fixed point and consider the corresponding Yoccoz puzzle.
%
%

Given a (geometric) puzzle piece $Y$ of depth $m$,
recall that $K_Y$ stands for $Y\cap K(f)$ and  $S_Y = \cl (K(f) \sm K_Y)$.
Let  $\YY$ stand for the space of paths $\de: [0,1] \to \BU^m \sm S_Y$
such that:
\begin{itemize}
  \item $\de(0)\in Y$,
  \item if $\de(t)\in K_Y$, then the restriction $\de|\, [0,t]$ is homotopic rel endpoints
       to a path contained in $Y$.\footnote{This condition can be replaced with a more restrictive one:
 After the first exit from $Y$, the path never intersect the Julia set $K(f)$
  (though it is allowed to return back to $Y$).}
\end{itemize}
Let $\BY$ be the space of paths $\de\in \YY$ modulo  homotopy through $\YY$ with $\de(1)$ fixed.
Define the projection $\pi_Y: \BY\ra \BU^m$ by  $[\de]\mapsto \de(1)$.
One can see that $\BY$ is a Riemann surface,
and $\pi_Y$ is an immersion such that
   $Y$ lifts to  a disk $\hat Y\subset \BY$ which is  homeomorphically projected onto $Y$.
 Thus, we can identify $\hat Y$ with $Y$; in particular,  $K_Y$ is embedded into $\BY$.

   The Riemann surface $\BY$ will be called the {\it pseudo-piece} (``$\psi$-piece'') associated with $Y$.

\ssk The $\psi$-pieces can also be defined in a different way.
Let us consider the topological annulus $A=\BU^m  \sm K(f)$ and its universal covering $\hat A$.
Let $Y_i$ be the components of $Y\sm K_Y$. There are finitely many of them,
and each $Y_i$ is simply connected. Hence they can be embedded into $\hat A$.
Select such an embedding $e_i: Y_i\ra \hat A_i$ where $\hat A_i$ stands for a copy of $\hat A$.
Then the $\psi$-piece is obtained by gluing the $A_i$ to $Y$ by means of $e_i$, i.e.,
$\BY= Y \sqcup_{e_i} \hat A_i.$

\begin{lem}
  The above two definitions of $\psi$-pieces are equivalent.
\end{lem}

\begin{proof}
  Let $\BY$ be a $\psi$-piece according to the first definition.
The puzzle piece $Y$ is embedded into $\BY$ by
associating to a point $y \in Y$ the constant path $\de(t)\equiv y$.

Let us realize the universal covering $\hat A_i\ra A$ as the space of paths
in $A$ that begin in $Y_i$ rel homotopy through such paths fixing the terminal endpoint.
(This realization is legitimate since $Y_i$ is simply connected.)
This provides us with an embedding $\phi_i: \hat A_i\ra \BY $

The embeddings $\phi_i$ have disjoint images.
Indeed, all points of $\di Y\cap K(f)$ are dividing
and thus belong to $S_Y$. Hence, if we take two paths $\de_1: [0,1]\ra \BU^m\sm K(F)$
and $\de_2: [0,1]\ra  \BU^m\sm K(f)  $ as above
representing points in $\hat A_i$ and $\hat A_j$ ($i\not= j$) with a common endpoint,
then they ``surround'' some piece of $S_Y$, and hence represent different points in $\BY$.

Moreover,  the image $\phi_i(\hat A_i)$ overlaps with  $Y$ by $Y_i$.
 Hence we obtain an embedding of $Y\sqcup_{e_i} \hat A_i$ into $\BY$.

Let us show that this embedding is surjective.
Take  a path $\de\in \tilde\YY$ representing some point of $\BY$,
and let $\tau\in [0,1]$ be the last parameter for which $\de(\tau)\in K_Y$.
Since the path $\de: [0,\tau]\ra \BU^m$ is trivial
(i.e., it can be pulled to $Y$ in $\BU^m\sm S_Y$ rel endpoints),
the restriction $\de: [\tau, 1]\ra \BU^m \sm S_Y$ (appropriately reparametrized)
represents the same point in $\BY$ as the original path.
Moreover, if $\tau\not=1$, we can replace it with an equivalent path
$\de: [\tau+\eps , 1]\ra \BU^m \sm S_Y$
which is disjoint from the Julia set $K(f)$.
As the latter path represents a point in some $\hat A_i$, we are done.
\end{proof}

\sss{Naturality}


\begin{lem}\label{psi-maps}
\begin{itemize}
\item [(i)]  Consider two puzzle pieces $Y$ and $Z$ such that the map $f: Y\ra Z$ is a branched covering
of degree $k$ (where $k=1$ or $k=2$ depending on whether $Y$ is off-critical or not).
Then there exists an induced map $\f: \BY \ra \BZ$
which  is a branched covering of the same degree $k$.

\item [(ii)] Given two puzzle pieces $Y\subset Z$, the inclusion $i: Y\ra Z$ extends to an immersion $\i: \BY\ra \BZ$.
\end{itemize}
\end{lem}

\begin{proof}
  Both properties follow easily from either definition of the $\psi$-pieces.
Let us, for instance, use the second definition.

(i)
Let $\depth Y=m$, $\depth Z=m-1$. Let us consider the degree $k$ branched covering
$$
    f: (\BU^m, Y, K_Y) \ra (\BU^{m-1}, Z, K_Z).
$$
The components $Y_i$ of $Y\sm K_Y$
are univalently mapped onto components $Z_{j(i)}$ of $Z\sm K_Z$, where the map
$j=j(i)$ is $k$-to-1.
This map extends to an isomorphism map $\hat A_i\ra \hat B_j $
of the corresponding universal coverings, which glue together into a branched covering
$\BY\ra \BZ$ of degree $k$.

(ii) Let $\depth Y=m$, $\depth Z=n$.
Let us  consider the immersion
$$
  i: (\BU^m, Y, K_Y) \ra (\BU^n, Z, K_Z).
$$
The components $Y_i$ are embedded by $i$ into some components $Z_{j(i)}$,
where the map $j=j(i)$ is surjective but not necessarily injective.
These embeddings extend to immersions $\hat A_i\ra \hat B_j$ that
glue together into an immersion $\BY\ra \BZ$.
\end{proof}

\sss{Moduli}

 Given two puzzle pieces $Z\Subset Y$, we let
$$
  \psimod (Y,Z) = \mod(\BY\sm K_Z).
$$

Lemma \ref{psi-maps} implies:

\begin{lem}\label{psi-mod}
\begin{itemize}
  \item [(i)]  Consider two pairs of puzzle pieces $(Y',Y)$ and $(Z',Z)$
such that the map $f: (Y',Y)\ra (Z',Z)$ is a branched covering
of degree $k$ (on both domains).
Then $$\psimod(Z',Z)=k\, \psimod(Y',Y).$$

\item [(ii)] Given a nest of three puzzle pieces $W\subset Z\subset  Y$,
we have $$\psimod(Z,W)\leq \psimod(Y,W).$$
\end{itemize}
\end{lem}

\subsubsection{Boundary of puzzle pieces}\label{boundary}
Let us mention in conclusion, that the ideal boundary of a pseudo-puzzle $\BY$
is tiled by (finitely many) arcs $\lambda_i\subset \di \hat A_i$
that cover the ideal boundary of $\BU^m$ (where $m=\depth Y$)
and arcs $\xi_i, \eta_i \subset \di \hat A_i$ mapped onto the Julia set $J(f)$.
The arc $\la_i$ meets each $\xi_i,\, \eta_i$ at a single boundary point
corresponding to a path $\de: [0,1)\mapsto A$ that wraps around $K(f)$ infinitely many times,
while $\eta_i$ meets $\xi_{i+1}$ at a vertex $v_i\in Y\cap K(f)$.
We say  that the arcs $\la_i$ form the {\it outer boundary } (or ``$O$-boundary'')
 $\di_O\BY$ of the puzzle piece $\BY$,
while the arcs $\xi_i$ and $\eta_i$
form its {\it $J$-boundary} $\di_J \BY$.
Given a vertex $v=v_i$ of a puzzle piece $Y$, let $\di^v\BY=\eta_i\cup \xi_{i+1}$
stand for the part of the $J$-boundary of $\BY$ attached to $v$.

Note that the immersion
constructed in Lemma \ref{psi-maps} extends continuously to the boundary of the
puzzle piece $\BY$.
(However, $\i(\di \BY)$ is not contained in $\di\BZ$, unless $Z=Y$.)
In what follows we will assume this extension without further comment.

A {\it multicurve} in some space $X$ is a continuous map $\gamma: \cup_{k=1}^l [s_k, t_k]\ra X$
parametrized by a finite union of disjoint intervals $[s_k, t_k]\subset \R$.%
\footnote{We allow that the boundary points of a  multicurve in a pseudo puzzle $\BY$ belong to  $\di_J\BY$.}
Note that  multicurves are ordered.
A multicurve in a puzzle piece $\BY$  is called {\it horizontal} if
$$
   \gamma(s_1)\in \di^{v_0} \BY, \quad  \gamma(t_k), \gamma(s_{k+1})\in \di^{v_k}\BY,\ k=1, \dots, l-1,
                  \quad \gamma(t_l)\in \di^{v_l} \BY
$$
for some vertices $v_k$ of $\BY$, $k=0, \dots, l$. We say that such a multicurve ``connects''
$\di^{v_0} \BY $ to $\di^{v_l} \BY$.
The following statement motivates introduction of  multicurves:

\begin{lem}\label{lifts}
Let $v$ and $w$ be two vertices of a geometric puzzle piece $Y\in\YY^m(l)$.
Then any curve $\gamma$ in $\BU^l$ connecting $S_Y^v$ to $S_Y^w$
contains a multicurve $\gamma'$  that lifts to a multicurve  $\gamma^*$ in $\BY$
connecting $\di_v\BY$ to $\di_w\BY$.
\end{lem}

Given two vertices $v$ and $w$ of $Y$, let $\GG_Y(v,w)$ stand for the family of horizontal multicurves in $\BY$
connecting $\di^v(\BY)$ to $\di^w(\BY)$.
Finally, let
$$
         \Bd_Y(v,w)= \LL(\GG_Y(v,w))
$$
stand for the extremal distance between the corresponding parts of $J$-boundary of $\BY$.

\begin{lem}\label{Bd}
  If  $f|Y$ is univalent,
then $\Bd_Y(v,w)=\Bd_{f(Y)}(fv, fw). $
\end{lem}

\comm{****************

\subsection{Pseudo-summary}\label{psi-summary}

   We can now recast the Summary of \S \ref{summary} in the pseudo-context.
We start with a $u\psi$-PL map $\bf$,
restrict it to a unicritical PL map $f$ and construct the associated dynasty of kingdoms.
We fix an arbitrary $m$, assume that the height of the dynasty is greater than $\log_2 m +5$,
and take some $n>\log_2 m +5$.
Then for any domain $\La =\La_k$, we obtain  the map $\Psi=\Psi_k$ \ref{3 domains}
satisfying properties (P1)-(P4). Applying to it the $\psi$-functor, we obtain:

\begin{itemize}

\item [(P1$'$)] An even-valued immersion ${\Bpsi}^{-1}: \BV^0 \ra \BUps$
 of degree bounded by $d^{2n+m}$;

\item[(P2$'$)] An even-valued immersion  ${\Bpsi}^{-1}: \BLa'\ra \BF$
  of degree bounded by  $d^5$;

\item [(P3$'$)] An immersion $\BUps\ra \BV^n$;

\item[(P4$'$)] An immersion $\hat\BE^{i+1} \ra  \BLa' $  which restricts to
   a homeomorphism $K_{E^{i+1}}\ra   K_\La$.

\end{itemize}

Moreover,

\begin{itemize}
\item[(P5$'$)] The images of the  immersions
 $\BLa_k'\ra \BV$ are disjoint from the Julia pieces $K_{\La_j}$,  $k\not=j$.

\item[(P6$'$)]
There are at least $m/2$ domains $\La_k$.

\end{itemize}
******************************************}

\section{Quasi-Additivity Law and Covering Lemma}\label{covering lemma sec}

Let us now formulate two analytic results
which will play a crucial role in what follows.
The first one appeares in \S 2.10.3 of \cite{covering lemma}:

\begin{qal}
    Fix some $\eta\in (0,1)$.
Let $\BV$ be a topological disk,  let  $K_i\Subset \BV $, $i=1,\dots, m$,
be pairwise disjoint full compact continua,
and let $\phi_i: \A(1,r_i)\ra \BV \sm   \cup K_j$
be holomorphic annuli such that each $\phi_i$ is
an  embedding of some proper collar of $\T$  to a proper collar of $\di K_i$.
Then there exists a $\de_0 > 0$ (depending on $\eta$ and $m$) such that: \\
If  for some $\de\in (0, \de_0)$,
$\mod (\BV, K_i) < \de$ while $\log r_i > 2\pi \eta\de$  for all $i$,
then
$$
  \mod(\BV,  \cup K_i ) < \frac{2 \eta^{-1} \de}{m}.
$$
\end{qal}

The next result appears in \S 3.1.5 of \cite{covering lemma}:

\begin{covlem}
Fix some  $\eta\in (0,1)$. Let us consider
two topological disks $\BU$ and $\BV$, two full continua $A'\subset \BU$ and $B'\subset \BV$,
and two compact subsets, $A\Subset A'$ and $B\Subset B'$, of topological type bounded by $T$.%
\footnote{In applications, $A$ and $B$ will be full continua, so $T=1$. }
\\
  Let $f:\BU \ra \BV$ be a branched covering of degree $D$ such that
$A'$ is a component of $f^{-1}(B')$, and $A$ is the union of some components of $f^{-1}(B)$.
Let $d=\deg(f : A'\ra B')$.\\
Let $B'$ be also embedded into another topological disk $\BB'$.
Assume $\BB'$ is immersed into
$\BV$ by a map $i$  in such a way that $i|\, B'=\id$, $i^{-1}(B')=B'$,
 and $i(\BB') \sm B'$ does not contain the critical values of $f$.  \\
Under the following  ``Collar Assumption'':
$$
   \mod(\BB', B) > \eta \mod(\BU,  A),
$$
if
$$ \mod(\BU, A) < \eps(\eta,T, D) $$
 then
$$
      \mod (\BV, B)  <  2 \eta^{-1} d^2 \mod(\BU,  A).
$$
\end{covlem}

\section{Improving the moduli}\label{moduli}

In this section $C$ will stand for the maximum of the constants in the Quasi-Additivity Law and the Covering Lemma.

\subsection{High type}

Let us begin with a simple estimate that compares moduli on consecutive odd levels of the Principal Nest:

\begin{lem}\label{simple}
For any odd $n$, we have:
 $$  \psimod(E^{n-3}, E^{n-2})\leq 4\, \psimod(E^{n-1}, E^n) $$
and
  $$ \psimod(Y^0, R) \leq 2^{\n+1}\, \psimod(E^0, E^1).$$
\end{lem}

\begin{proof}
By Lemma \ref{psi-mod},
$$
   \psimod(\hat E^{n-1}, E^{n-1}) = 2\, \psimod(\tilde E^n, E^n) \leq 2\, \psimod(E^{n-1}, E^n)
$$
and
$$
   \psimod(\hat E^{n-1}, E^{n-1}) \geq \psimod(\tilde E^{n-1}, E^{n-1} )
$$
$$
=\frac{1}{2}\, \psimod(\hat E^{n-2}, E^{n-2})=
                                     \frac{1}{2}\,  \psimod(E^{n-3}, E^{n-2}),
$$
and the first estimate follows.

The second estimate is similar. The puzzle piece $E^1\equiv Y^m$
is mapped with degree 2 onto $E^0$, and this map admits degree 2 extension
$\tilde E^1\ra Y^{q\n}\equiv \hat E^0$, where $\tilde E^1 = Y^{m-1}$.
Then $E^0$ is mapped onto $Z^1_\kappa$ by degree 2 map $f^{q\n}$.
This map admits degree $2^{\n}$ extension $Y^{\q\n}\ra Y^0$.
It follows that
$$
   \psimod(Y^0, R)\leq \psimod(Y^0, Z^1_\kappa)
$$
$$
  = 2^{\n+1}\,\psimod(Y^{m-1}, E^1)\leq 2^{\n+1}\, \psimod(E^0, E^1).
$$
\end{proof}

The following  lemma tells us that if some principal  modulus is very small then it
should be even smaller on some preceding level of the Principal Nest:

\begin{lem}\label{growth-lem}
There exist absolute  $N\in \N$ and $\eps>0$  such that:
If on some odd level $n\geq N$, $  \psimod  (E^{n-1}, E^n ) <\eps $,
then on some previous odd level $n-s\in [n-N, n-1]$ we have:
\begin{equation}\label{decrease}
    \psimod  (E^{n-s-1},  E^{n-s} ) <  \frac{1}{2}\, \psimod( E^{n-1},  E^n ).
\end{equation}
\end{lem}

\begin{proof}
Let us fix some integer $m >  C^3 2^{28}$.
Let $N$ be the smallest odd integer that is bigger than $\log_2 m + 5$.
Take any odd level $n\geq N$. For each $k$,
let us consider the associated 3-domain branched covering $\Psi_k$ (\ref{3 domains})
$$
   \Psi_k: (\Upsilon_k, \De_k, E^n) \ra (E^{n-N}, \La_k', \La_k).
$$
Let us consider two cases:

\msk {\it Case 1.} Assume that for some domain $\La_k$,
$$
  \psimod (\La_k', \La_k)\leq  \frac{1}{4}\, \psimod (E^{n-1},  E^n).
$$
         By Property (P4) and Lemma \ref{psi-mod},
$\psimod(\La_k' ,  \La_k) = \psimod (\hat E^i, E^i)$.
If $i$ is odd then $\hat E^i= E^{i-1}$,
 and   we obtain the desired estimate with $s=n-i\in [1, 5]$:
$$
  \psimod (E^{n-1}, E^n) \geq 4\, \psimod ( E^{i-1}, E^i).
$$
If $i$ is even, then
$$
     \psimod (\hat E^i, E^i) \geq \psimod(\tilde E^i, E^i) = \frac{1}{2} \psimod ( \hat E^{i-1}, E^{i-1})=
                              \frac{1}{2} \psimod ( E^{i-2}, E^{i-1}),
$$
and we conclude that
$$
     \psimod (E^{n-1}, E^n) \geq 2\, \psimod ( E^{i-2}, E^{i-1}).
$$

\msk{\it Case 2.} Assume that for all $\La_k$,
\begin{equation}\label{1/2}
  \psimod (\La_k', \La_k) \geq \frac{1}{4} \psimod (E^{n-1},  E^n)\geq \frac{1}{4} \psimod(\Upsilon_k, E^n)
\end{equation}
(where the second estimate follows from the inclusion $\Upsilon_k\subset E^{n-1}$).
By Lemma \ref{psi-maps}, there exists a natural covering map
$$
   \boldsymbol{\Psi}_k: (\boldsymbol{\Upsilon}_k, K_{\De_k}, K_{E^n} )\ra (\BE^{n-N-1}, K_{\La_k'}, K_{\La_k}),
$$
and  a natural immersion $ i: \boldsymbol{\La}_k'  \ra \BE^{n-N} $.
Note that $i(\boldsymbol{\La}_k')\sm K_{\La_k'}$ does not contain the critical values of $\boldsymbol{\Psi}_k$,
since the latter are contained in the Julia set $K(f)$.
Moreover, equation (\ref{1/2}) provides us with the Collar Assumption
that allows us to apply the Covering Lemma to the map $\boldsymbol{\Psi}_k$.
If $\eps$ is sufficiently small, it yields:
\begin{equation}\label{outer mod}
  \psimod (E^{n-N-1},  \La_k) \leq C 2^{12} \psimod(\Upsilon_k,  E^n) \leq C\, 2^{12}\, \psimod(E^{n-1}, E^n).
\end{equation}
Estimates (\ref{1/2}) and (\ref{outer mod}) show that the Quasi-Additivity Law is applicable
to the family of islands $K_{\La_k}$ in $\BE^{n-N-1}$
with $\eta^{-1}= C 2^{14}$.
 Since there are at least $m/2$ domains $\La_k\subset \La_k'\subset E^{n-N}$,
it implies:
$$
      \psimod(E^{n-N-1}, E^{n-N}) \leq \frac{C^3 2^{27} \, \psimod(E^{n-1}, E^n )}  {m}  < \frac{1}{2}\, \psimod(E^{n-1}, E^n),
$$
and we are done.
\end{proof}

\begin{lem}\label{exp decay}
   There exist absolute constants $C>0$, $\rho\in (0,1)$ and $\eps>0$
such that if for some odd $n$,  $\psimod(E^{n-1}, E^n)<\eps$,
then
$$
   \psimod(E^0, E^1)< C \rho^n \psimod(E^{n-1}, E^n).
$$
and
$$
   \psimod(Y^0, R)\leq  C 2^\n \rho^n \psimod(E^{n-1}, E^n).
$$
\end{lem}

\begin{proof}
  By Lemma \ref{growth-lem},
there exists an odd level $l<N$ such that
$$
   \psimod(E^{l-1}, E^l) \leq \left(\frac{1}{2}\right)^{[n/N]},
$$
which together with Lemma \ref{simple} implies the desired estimates.
\end{proof}

\subsection{Frequent $R$-returns}

 Let us consider the map
\begin{equation}\label{l}
    f^l= f^{q\n+1} \circ g_1\circ\dots  \circ g_{\chi-1}\, :\, E^{\chi-1}\ra Y^0
\end{equation}
and the trajectory $\OO= \{\KK_i\}_{i=l}^{l+p-1}$
of the little Julia set $\KK$. Let $i_1, i_2,\dots$ be the moments in $\OO$
for which $\KK_i\subset  R$.

\begin{lem}\label{frequent returns}
Let $\rho>0$,  $\bar\chi\in \N$.
Take some integer $m\geq C^3 2^{30}/\rho$, and let $\underline p= m^2\q\n$.
Assume that the little Julia set frequently visits $R$:
\begin{equation}\label{frequent visits}
    i_{k+1}- i_k \leq m\q\n ,\quad k=1,2,\dots, m.
\end{equation}
If $\chi\leq \bar\chi$ while $p\geq \underline p$,
then
$$
  \psimod (Y^0, R) \leq \rho \mod(\BE^{\chi-1}, \KK),
$$
provided $\mod(\BE^{\chi-1}, \KK)< \eps(\n; \bar\chi, \rho)$.
\end{lem}

\begin{proof}
The map (\ref{l})
has degree  $2^\chi\leq 2^{\bar\chi}$. By Lemma 2.9 of \cite{KL-high}, \note{formulate separately}
 $\deg (f^l |\, E^\chi)\leq 32$,  and hence $l\leq 5p$.

 By (\ref{frequent visits}),   $i_m-i_1 < m^2 \q \n \leq p $,
so the moments $i_k$ are pairwise non-congruent mod $p$.
Hence the little Julia sets $\KK_{i_1}, \dots, \KK_{i_m}$
are all distinct.

Since $\OO$ has length $p$, there is only one critical Julia set in $\OO$.
Hence $\deg(f^{i_k}: \KK\ra \KK_{i_k}) $ is at most 64, so that $i_k\leq 6p$,
$k=1,\dots, m$.

On the other hand, $\KK_{i_k}$ is contained  in  a puzzle piece in $R$
which is mapped under $f^{i_{k+1}-i_k}$ onto $Y^0$ with degree at most $2^{m\n}$. \note{lemma}
It follows by the Telescope Lemma \ref{telescope} that there is a puzzle piece $\Upsilon_k\subset E^{\chi-1}$
which is mapped under $f^{i_k}$ onto $Y^0$ with degree at most $2^{\bar\chi+km\n}\leq 2^{\bar\chi + m^2\n}\equiv D$.

We would like to apply the Covering Lemma to the corresponding map
$$
   \f^{i_k}: (\boldsymbol{\Upsilon}_k, \KK)\ra  (\BY^0, \KK_{i_k})
$$
of degree at most $D$.
To this end we need collars around $\KK_{i_k}$.
Let $\Om$ be the critical pullback of  $E^\chi$ under $f^{6p}$.
 Then we let $\La_k' = f^{i_k}(\Om))$.
Since the moments $i_k$ are pairwise non-congruent mod $p$ and $i_k\leq 6p$,
 the puzzle pieces $\La_k'$ are contained in different domains of the orbit
$f^t(E^\chi)$, $t=0,1\dots, p-1$. Hence they are pairwise disjoint.
Moreover, by Lemma \ref{psi-mod},
\begin{equation}\label{collars again}
  \mod(\BLa_k', \KK_{i_k})\geq  \mod(\boldsymbol{\Om}, \KK)= \frac{1}{4} \mod(\BE^{\chi-1}, \KK) \geq  \frac{1}{4} \mod(\BUps_k, \KK).
\end{equation}
This provides us with the desired Collar Assumption.
  By the Covering Lemma,
$$
  \mod(\BY^0,\KK_{i_k})\leq C 2^{14} \mod(\BUps_k, \KK)\leq C 2^{14} \mod(\BE^{\chi-1}, \KK) .
$$
The last two estimates show that
the Quasi-Additivity Law is applicable to the family of islands $\KK_{i_k}$ in $\BY^0$
(with $\eta^{-1}= C2^{16}$):
$$
   \psimod (Y^0, R) \leq \frac{C^3 2^{30} \mod(\BE^{\chi-1}, \KK)}{m}\leq \rho \mod(\BE^{\chi-1}, \KK),
$$
provided $\mod(\BE^{\chi-1}, \KK)< \eps(D) = \eps(\n; \bar\chi, \rho)$
and we are done.
\end{proof}

\subsection{Many consecutive returns to $L$}\label{Robespier}

Here the set-up is the same as in the previous section,
but we will assume that there is a gap in returns of the little Julia sets to $R$:

\begin{lem}\label{Robespier lem}
Let $\rho$,  $\bar\chi$, $m$, and $\underline p$ be as in Lemma \ref{frequent returns}.
Assume there is $k\leq m$ such that
\begin{equation}\label{unfrequent visits}
    i_{k+1}- i_k > m\q\n .
\end{equation}
If $\chi\leq \bar\chi$ while $p\geq \underline p$,
then
$$
  \psimod (Z^0, L) \leq \rho \mod(\BE^{\chi-1}, \KK),
$$
provided $\mod(\BE^{\chi-1}, \KK)< \eps(\n; \q, \bar\chi)$.
\end{lem}

\begin{proof}
Under our assumption (\ref{unfrequent visits})
the Julia set returns frequently to $L$:
$$
        \KK_{i_k+ j \q\n} \subset L, \quad j=1,\dots, m.
$$
Let $P_j \ni f^{i_k}z $ be the pullback of $Z^0$ under $f^{j \q\n}$.
By Lemma \ref{qn-pullbacks} and the Telescope Lemma,
$P_j\subset Y^0$.

Let $\Ups_j$ be the further pullback of $P_j$ under $f^{i_k}$,
and let
$$
    \Psi_j= f^{i_k+j\q\n}: \Ups_j\ra Z^0.
$$
Then $\Ups_j\subset E^{\chi-1}$ and $\deg \Psi_j\leq 2^{\bar\chi+\underline p}$.

The rest of the argument is the same as for Lemma \ref{frequent returns}:
the Covering Lemma implies that for $j=1,\dots, m$,
$$
    \mod (\BZ_0, \KK_{i_k+j\q\n}) \leq C 2^{14} \mod(\BUps_k, \KK)\leq  C 2^{14} \mod(\BE^{\chi-1}, \KK) ,
$$
and by the Quasi-Additivity Law,
$$
  \psimod(Z^0, L) \leq \frac{C^3 2^{30} \mod(\BE^{\chi-1}, \KK)}{m}\leq \rho \mod(\BE^{\chi-1}, \KK).
$$
\end{proof}

Note that by symmetry, $\psimod(Z^0, L) = \psimod(Y^0, R)$.
Putting together Lemmas \ref{exp decay}, \ref{frequent returns} and \ref{Robespier lem},
we obtain:

\begin{cor}\label{R-L dist}
For any parameters $\q,\n$ of a decoration and any $\rho>0$,
 there exists  $\underline p \in \N$ and $\eps>0$ such that
$$
  \psimod (Y^0, R) \leq \rho \mod(\BE^{\chi-1}, \KK),
$$
provided  $p\geq \underline p$ and $\mod(\BE^{\chi-1}, \KK)< \eps$.
\end{cor}


\subsection{Comparison of $\psimod(Y^0, R)$ with $\Bd_{Y^1}(\alpha, \alpha')$}

  Let  $$\mu := \min (\psimod(U, K), 1/2) .$$

\begin{lem}\label{observe}
    $$\Bd_{Y^1}(\alpha, \alpha')\leq \psimod (Y^0, R) \ominus \frac{1}{2} \mu.$$
\end{lem}

\begin{proof}
  The boundary of $\BY^0$ consists of two parts (see the end of \S \ref{psi-puzzle}):
the $J$-boundary $\di_J\BY^0 = \xi\cup\eta$ attached to $\alpha$ and the outer  arc $\la=\di_O\BY^0$
that covers the ideal boundary of $\BU$.
Let $\GG^h$ stand for the family of curves in the annulus $\BY^0\sm K_R$  connecting $K_R$ to the $J$-boundary,
while $\GG^v$ stand for the family of curves in the same annulus  connecting $K_R$ to $\la$.
By the Parallel Law,
$$
   \LL(\GG^h) \leq \psimod (Y^0, R) \ominus \LL(\GG^v).
$$
Let $\Pi$ stand for the rectangle uniformizing  $\BY^0\sm K_{Y^0}$ whose horizontal sides
correspond to  $K_{Y^0}$ and $\la$,
and vertical sides correspond to $\xi$ and $\eta$.
We let $\om$ be the horizontal side of $\Pi$ corresponding to $K_{Y^0}$.
Since any curve of the family $\GG^v$
overflows some curve connecting $K_{Y^0}$ to $\la$ in $\BY^0\sm K_{Y^0}$
(and thus representing a vertical curve in $\Pi$), we have:
$$
\LL (\GG^v) \geq \mod \Pi.
$$
But by definition of the pseudo-puzzle,
the domain $\BY^0\sm K_{Y^0}$  covers the annulus $\BU\sm K$
extending to an embedding on $K_{Y^0}$. Let us uniformize $\BU\sm K$ by a
round annulus $\A$. It follows that the rectangle $\Pi$ covers $\A$ in such a way that
$\om\subset \di\Pi $ is embedded into $\di \A$.
By Lemma~\ref{C} from the Appendix,
$$
   \mod \Pi \geq \frac{\mu}{2}.
$$
Putting the above three estimates together, we obtain:
\begin{equation}\label{pl}
  \LL(\GG^h) \leq \psimod (Y^0, R) \ominus \frac{\mu}{2}.
\end{equation}

On the other hand,
let us consider the family $\HH$ of horizontal curves in the puzzle piece $\BY^1$
connecting $\di_\alpha\BY^1$ to $\di_{\alpha'} \BY^1$.
Let $\phi: \BY^1\ra \BY^0$ be the natural immersion.
Under $\phi$,
the boundary $\di_\alpha \BY^1$ is mapped homeomorphically onto $\di_\alpha \BY^0$.
It follows that any curve $\gamma$ of $\GG^h$ contains an arc
that can be lifted by $\phi$ to some curve of $\HH$.
Indeed, orient $\gamma$ so that it begins on $\di_\alpha\BY^0$.
Then a maximal lift of $\gamma$ that begins on $\di_\alpha\BY^1$
must end on $\di_{\alpha'} \BY^1$.

By Corollary \ref{two families}, $ \LL(\HH)\leq \LL(\GG^h)$.
Together with (\ref{pl}), this yields the desired inequality.
\end{proof}

\subsection{Skipping over}
  In this section we will show that not many curves can skip some piece of the Julia set.

Let $Y\in \YY^m(l)$ be a geometric puzzle piece of bidepth $(m,l)$,
and let $A$ be a component of $Y\sm K_Y$. Let $\di_J \hat A = \di_J \BY \cap \hat A $.
Recall that it consists of two components.
Let $\CC_A$ stand for the family of curves in $\hat A$
connecting  different  components of $\di_J \hat A$,
and let
$$
  \Bd_A = \LL(\CC_A).
$$

Let $\displaystyle{\frac{1}{2\pi}\log r= 2^{-(m+\q)}\mu}$.
Then the annulus $\BU^{m+\q}\sm K$ can be uniformized by the round annulus
$\A(1,r)$, and under this uniformization,
the set $K_Y$ gets represented on the unit circle $\T$ as the union of arcs $\om_i$
of length
\begin{equation}\label{om}
    |\om_i|  \geq 2\pi\cdot 2^{-(\q+m)} .
\end{equation}
Indeed, the covering map $\Bf^{m+\q}: \BU^{m+\q} \sm K\ra \BU\sm K$
is turned into $z\mapsto z^{2^{m+\q}}$ under the above uniformization of $\BU^{m+\q} \sm K$
and the uniformization of $\BU\sm K$ by $\A(1, e^{2\pi\mu})$ (appropriately normalized).
Since under this map, every arc $\om_i$ covers the whole circle,
the length of $\om_i$  is at least $2\pi$ times its combinatorial length (\ref{comb length}).

\begin{lem}\label{m-m}
  Let  $Y\in \YY^m(m+\q)$ be a geometric puzzle piece of bidepth $(m,m+\q)$,
and let $A$ be a component of $Y\sm K_Y$. Then
$$
   \Bd_A \geq \frac{1}{\mu}.
$$
\end{lem}

\begin{proof}
%
We can uniformize $\hat A$ by the horizontal  strip $\Pi=\Pi(2^{-(m+\q)})\mu$
 in such a way that the upper boundary
of $\Pi$ covers the $O$-boundary of $\BU^{m+\q}$, and the group of deck transformations is generated by
the translation $z\mapsto z+1$. By (\ref{om}),
the Julia set $K_Y\subset \di A$ is represented
as an interval $I$ on $\R$ of length at least $2^{-(\q+m)}$.

Let us view $\Pi$ as a quadrilateral with horizontal sides $I$ and the top of $\Pi$.
Then
$$
   \LL(\CC_A)= \frac{1}{\mod \Pi} \geq \frac{1}{\mu },
$$
where the last estimate comes from the simple right-hand side estimate  of Lemma \ref{Pi},
and we are done.
%
\end{proof}

\begin{lem}\label{m-0}
  Let $Y\in \YY^m(0)$ be a geometric puzzle piece of bidepth $(m,0)$,
and let  $A$ be a connected component of $Y\sm K_Y$.
Then
$$
   \Bd_A \geq \frac{\mu}{2^{m+\q+2}}.
$$
\end{lem}

\begin{proof}
Let $A(m+\q)$ be the component of $Y(m+\q)\sm K_Y$ contained in $A$,
and let $\phi: \hat A(m+\q)\ra \hat A$ be the natural immersion.
It extends to the identity on $K_Y\cap \di \hat A(m+\q)$.

 Let us realize   $\hat A$ as the strip $\Pi(\mu)$ that covers $\BU \sm K$,
with the group of deck transformations generated by $z\mapsto z+1$.
Let us consider the interval $I \subset \R$  representing $K_Y\cap \di A$,
and let $J$ be the left-adjacent interval of length 1.

Let us orient the curves $\gamma\in \CC_A$ so that they begin on the left-hand side of $I$.
Then any curve $\gamma\in \CC_A$ contains the maximal initial arc $\gamma'$ that can be lifted
by $\phi$ to a curve $\gamma^*$ in  $\hat A(m+\q)$.
Accordingly,  we can split the family of curves $\gamma\in \CC_A$ into three subfamilies:
\begin{itemize}
  \item $\HH_1$ consists of the curves $\gamma$ such that $\gamma'=\gamma$;
   then  $\gamma^*\in \CC_{A(m+\q)}$;

  \item $\HH_2$ consist of the curves that begin in $J$ and whose
lift $\gamma^*$ terminates on the $O$-boundary of $\hat A(m+\q)$;

   \item $\HH_3$ consists of the curves that begin on  the left-hand side of $J$.
\end{itemize}

Let us estimate the extremal length of each of these families.

Since $\HH_1 = \phi (\CC_{A(m+\q)})$,
$$
  \LL(\HH_1)\geq \LL(\CC_{A(m+\q)})\geq \frac{1}{\mu},
$$
where the first estimate follows from Lemma \ref{increase}, and the second follows from
Lemma \ref{m-m}.

Let $\TT$ be the family of curves in $\hat A(m+\q)$ that begin on $J$ and end on the
$O$-boundary of $\hat A(m+\q)$. By Corollary \ref{two families} and Lemma \ref{Pi},
$$
   \LL(\HH_2)\geq \LL(\TT)\geq \frac{\mu}{2^{m+\q+1}}.
$$

To estimate the extremal length of $\HH_3$, endow the rectangle
$Q=J \times [0,\mu]\subset \Pi(\mu)$ with the Euclidean metric $\la$.
Since  any curve $\gamma\in \HH_3$ horizontally overflows $Q$,
it has $\la$-length at least 1. Hence
$$
   \LL(\HH_3)\geq \frac{1}{\area Q}= \frac{1}{\mu}.
$$

Incorporating the last three estimates into the Parallel Law, we obtain the desired:
$$
    \LL(\CC_Y^\eps) \geq \frac{1}{\mu + 2^{m+\q+1}\mu^{-1} +\mu}\geq \frac{\mu}{2^{m+\q+ 2}}.
$$
%
\end{proof}

 Let us consider two vertices,  $v$ and $w$, of a geometric puzzle piece $Y$.
   Let $Z\subset Y$ be a  puzzle piece of depth $m$ that separates $v$ from $w$.
We say that a multicurve in $\BY$ connecting $\di^v \BY$ to $\di^w\BY $ {\it skips over} $K_Z$
if one of its components does not cross $K_Z$.

\begin{cor}\label{separating pieces}
Under the above circumstances,
let $\TT$ be be the family of multicurves in $\BY$ connecting $\di^v \BY$ to $\di^w\BY $
that skip over   $K_Z$.
Then
$$
  \LL(\TT)\geq \frac{\mu}{2^{2 m+\q+2} }.
$$
\end{cor}

\begin{proof}
  The piece $Z$ has at most $2^m$ components $A_i$ of $Z\sm K_Z$.
If a multicurve $\gamma$ in $\BY$ that skips over $K_Z$
then it contains an arc $\gamma'$ that lifts to
a curve $\gamma^*$ in some family $\GG_{A_j}$.
Let $\CC_j$ be the corresponding subfamily of $\TT$.
By Lemma \ref{m-0} (together with Corollary \ref{two families}),
$$
    \LL(\CC_j) \geq \Bd_{A_j} \geq \frac{\mu}{2^{m+\q+2}}.
$$
The Parallel Law concludes the proof.
\end{proof}

 Let us now  consider the puzzle piece $P=Y^{(\n-1)\q+1}$,
together with the corresponding pseudo-piece $\BP$,
and the family of puzzle pieces $Q^v\subset P$
from   Lemma \ref{separating family}.
Recall that $T^{vw}= K_P\sm (Q^v\cup Q^w)$.
Given two vertices $v$ and $w$ of $P$,
let $\hat\GG_P^{vw}$ stand for the family of multicurves in $\BP$ connecting $\di_v \BP$ to
$\di_w\BP$ that {\it do not  skip over} $T_{vw}$.
By Corollary \ref{separating pieces},
\begin{equation}\label{T^vw}
  \LL(\GG\sm \hat\GG_P^{vw}) \geq C^{-1} \mu,
\end{equation}
where here and below, $C$ stands for a  constant that depend only on $\q$ and $\n$.

\subsection{Separation between $L$ and $R$}

   In this section we will show that   the modulus $\Bd_{Y^1}(\alpha, \alpha')$
that measures  the extremal distance between $L$ and $R$ is comparable
with $\mu$.

Let $Y$ be a geometric puzzle piece.  For  two vertices $v$ and $w$ of $Y$,
 we let
$$
 \BW_Y(v,w) =\WW(\GG_Y^{vw}).
$$
We define the {\it pseudo-conductance} of $Y$ as
$$
    \BW_Y =\sup_{v,w} \BW_Y(v,w),
$$
where the supremum is taken over all pairs of the vertices of $Y$.

\begin{lem}
  For the puzzle piece $P=Y^{(\n-1)\q+1}$ we have:
   $$   W_P\leq \frac{C} {\mu}. $$
\end{lem}

\begin{proof}
    Along with the above conductance of $P$, let us consider
$$
    \hat\BW_P(v,w) =\WW(\hat\GG_P^{vw});\quad \hat\BW_P =\sup_{v,w} \hat\BW_P(v,w).
 $$
By (\ref{T^vw}),
\begin{equation}\label{hat}
    \BW_P \leq \hat\BW_P + \frac{C}{\mu}.
\end{equation}

Take a pair of vertices, $v$ and $w$.
Let $Q^v\cap T^{vw}=\{v'\}$ and  $Q^w\cap T^{vw}= \{w'\}$.
Recall that depth of the puzzle pieces $Q^v$ and $Q^w$ is equal to $r=(2\n-1)\q+1$,
and so depends only on $\q$ and $\n$. Let $\EE^r$ be the lift of the equipotential of level $r$
to $\BP$.

For any horizontal multicurve $\gamma\in \hat\GG_P^{vw}$,
one of the following two  possibilities can occur:

\ssk\nin $\bullet$
  $\gamma$ crosses the equipotential $\EE^r$, and hence it contains an arc $\gamma'$ connecting $\EE^r$
to $T^{vw}$; By Lemma \ref{C} and the Parallel Law,
the width of this family of curves is bounded by $2^{r+\n} / \mu$ (here $2^{\n-1}$
is a bound on the number of connected components of $P\sm K_P$);

\ssk\nin $\bullet$
 contains two disjoint multicurves,
$\de^v$ and $\de^w$, that do not cross $\EE^r$ and such that $\de^v$ connects $\di^v\BP$ to $T^{vw}$,
while $\de^w$ connects $T^{vw}$ to $\di^w\BP$. Then $\de^v$ contains a multicurve that can
be lifted to a horizontal multicurve in $\BQ^v$ connecting $\di^v\BQ^v$ to $\di^{v'}\BQ^v$,
 an similarly for  $\de^w$.

\ssk
By the Series and Parallel Laws,
$$
     \hat\BW_P(v,w) \leq \BW_{Q^v}(v, v') \oplus \BW_{Q^w}(w,w')+ \frac{2^{r+\n}}{\mu}
           \leq \BW_{Q^v} \oplus \BW_{Q^w}+\frac{2^{r+\n}}{\mu}.
$$
But $\BW(Q^v)=\BW(Q^w)= \BW(P)$  since $Q^v$ and $Q^w$ are univalent pullbacks of $P$.
Hence
$$
   \BW_{Q^v}\oplus \BW_{Q^w}  \leq \frac{1}{2} \BW_P.
$$
Putting the last two estimates together and taking the supremum over all pairs of vertices $(v,w)$ of $P$,
we conclude that
$$
   \hat\BW_P\leq \frac{1}{2}\BW_P + \frac{2^{r+\n}}{\mu}.
$$
 Together with  (\ref{hat}) it yields:
$$
   \BW_P\leq \frac{1}{2} \BW_P + \frac{C}{\mu},
$$
and the conclusion follows.
\end{proof}

\begin{prop}\label{separation prop}
  $\Bd_{Y^1}(\alpha, \alpha') \geq C^{-1}\mu. $
\end{prop}

\begin{proof}
   Since the map $f^{(\n-1)\q}: \BP\ra \BY^1$ is a branched covering
that maps $\di_0\BP$ to  $\di_0\BY^1$, any curve
$\gamma\in \GG_{Y^1}^{\alpha \alpha'}$ can be lifted to a curve
$\gamma^*\in \cup\GG_P^{vw}$, where the union is taken over all pairs of vertices of $P$.
Hence
$$
 \LL (\GG_{Y^1}^{\alpha \alpha'}) \geq \bigoplus_{v,w} \LL(\GG_P^{vw})\geq \frac{1}{N\BW_P},
$$
where $N$ is the number of pairs $(v,w)$.
The conclusion follows.
\end{proof}

 Lemma \ref{observe} and Proposition~\ref{separation prop} imply:

\begin{cor}\label{mu vs Y-R}
   $\displaystyle{\frac{\mu}{2} \leq C\, \psimod (Y^0, R).}$
\end{cor}

\begin{cor}\label{slow change}
Let $f: (\BU, K)\ra (\BU, K)$ be a renormalizable $\psi$-quadratic-like map
with decoration  parameters $(\q, \n)$,
and let $f'=f^p : (\BU', K')\ra (\BU',K')$ be its first renormalization.
Then
$$
      \min\{ \mod (\BU, K), 1/2\} \leq C\,  \mod(\BU', K'),
$$
where $C=C(\q,\n)$.
\end{cor}

\begin{proof}
  This follows from Lemma \ref{simple} and Corollary \ref{mu vs Y-R}
by noticing that $(\BE^{\chi-1}, \KK) = (\BU', K')$.
\end{proof}

\subsection{Conclusion}\label{conclusion}

Everything is now prepared for the main results.
 Corollary \ref{R-L dist} and Corollary \ref{mu vs Y-R} imply:

\begin{thm}[Improving of the moduli: bounded decoration parameters]\label{mod improve}
For any parameters $\bar\q,\bar\n$ and any $\rho>0$,
 there exist  ${\underline p} \in \N$ and $\eps>0$ with the following property.
Let $f: (\BU, K)\ra (\BU, K)$ be a renormalizable $\psi$-quadratic-like map
with  decoration parameters $(\q,\n)\leq (\bar\q, \bar\n)$,
and let $f'=f^p : (\BU', K')\ra (\BU',K')$ be its first renormalization.
Then
$$
 \{ p\geq \underline p\ \ \mathrm{and}\ \mod(\BU', K')< \eps \}  \imply \mod (\BU, K) < \rho \mod(\BU', K').
$$
\end{thm}


\begin{rem}  The logic of this theorem can be adjusted so that it would sound more
like an ``improvement in the future'' rather than ``worsening in the past'':

{\it  For any parameters $\bar\q, \bar\n$ of a Misuirewicz limb,
 there exists  $\underline p \in \N$ and $\eps>0$ such that
$$
    \mod(\BU', K') \geq   2\,  \mod (\BU, K)
$$
provided  $p\geq \underline p$ and $\mod (\BU, K) < \eps/2$.}
\end{rem}

\msk
  Theorem \ref{mod improve},
together with Lemma \ref{ql extension},
implies Theorem \ref{high periods thm} from the Introduction.

To derive the Main Theorem, we will combine Theorem \ref{mod improve}
with the following result  (Theorem 9.1 from \cite{K}):

\begin{thm}[Improving of the moduli: bounded period]\label{improving mod-2}
For any $\rho\in (0,1)$, there exists $\underline p=\underline p (\rho)$
such that for any $\bar p\geq \underline p$, there exists $\eps=\eps(\bar p)>0 $ with the following property.
Let $f: (\BU,K)\to (\BU,K)$ be primitively renormalizable $\psi$-quadratic-like map,
and let $f'=f^p = (\BU',K') \ra (\BU',K')$ be the corresponding renormalization.
Then
$$
   \{\underline p\leq p\leq \bar p\ \ \mathrm{and}\ \mod(\BU'\sm K') < \eps\}   \imply \mod(\BU\sm K) < \rho \mod(\BU'\sm K').
$$
\end{thm}

\begin{rem}
  Unlike Theorem \ref{mod improve}, in Theorem \ref{improving mod-2} the map
$f'$  {\it is not necessarily the first renormalization of} $f$.
On the other hand, in Theorem  \ref{improving mod-2},
the scale $\eps$ {\it depends on the upper bound} $\bar p$,
while in Theorem \ref{mod improve} it does not.
\end{rem}

We say that an infinitely renormalizable $\psi$-ql map $f$
belongs to the {\it decoration class} $(\bar\q, \bar\n)$ if the decoration parameters $(\bar\q_n, \bar\n_n)$
of the renormalizations $R^n f$ are all bounded by $(\bar\q, \bar\n)$.

Let us now put the above two theorems together:

\begin{cor}\label{improving-3}
 For any $(\bar\q, \bar\n)$,   there exist an $\eps>0$ and $l\in \N$ with the following property.
For any  infinitely renormalizable $\psi$-ql map $f$  of decoration class $(\bar\q, \bar \n)$
with renormalizations $R^n f: (\BU_n, K_n) \ra (\BU_n, K_n)$,
if $\mod(\BU_n\sm K_n)< \eps$, $n\geq l$,  then $\mod (\BU_{n-l}\sm K_{n-l})<\mod(\BU_n\sm K_n) /2$.
\end{cor}

\begin{proof}
  Given an infinitely renormalizable $\psi$-ql map $f$ with renormalizations
$R^n f: (\BU_n, K_n)\ra (\BU_n\sm K_n)$, we let $\mu_n(f)= \mod (\BU_n, K_n).$
Arguing by contradiction, we find a sequence $f_i$ of infinitely renormalizable $\psi$-ql maps
of decoration class $(\bar\q, \bar\n)$ and sequences $\eps_i\to 0$ and $n(i) \to \infty$
such that:

\ssk\nin P1:
   $\mu_{n(i)}(f_i) <\eps_i$;

\ssk\nin P2:
  $ \mu_{n(i)}(f_i) < 2 \mu_k (f_i)$,  $k=0,1,\dots,  n(i)-1 $.
\ssk

Let $R^n f_i $ be the renormalization of $R^{n-1} f_i$ with period $p_n(f_i)$.
Applying  then the diagonal process,  we can also assume the following property:

\ssk\nin P3:
   $p_{n(i)-s}(f_i)\to \pi_s\in \N \cup \{\infty\}$ for  $s=0,1, \dots$.
\ssk

We  let $\bar s \in \Z_{\geq 0} \cup\{\infty\}$ be the first moment for which $\pi_s=\infty$
(with understanding that $\bar s=\infty$ if such a moment does not exist).

Let us consider two cases:

\ssk
{\it Case 1}: $\bar s<\infty$.
Applying consecutively Corollary \ref{slow change}, we conclude that for sufficiently big $i$,
$$
   \mu_{n(i)-s}(f_i) \leq C^s \mu_{n(i)}(f_i), \quad s\leq \bar s.
$$
Let $\rho\in (0, 1/2 C^{\bar s})$.
By Theorem \ref{mod improve}, for all sufficiently big $i$,
$$
   \mu_{n(i)-\bar s-1}(f_i) \leq \rho\,  \mu_{n(i)-\bar s}(f_i).
$$
 Putting the last two estimates together, we conclude that for all sufficiently big $i$,
$$
   \mu_{n(i)-\bar s-1}(f_i) < \frac{1}{2} \mu_{n(i)}(f_i),
$$
  contradicting assumption (P2).

\ssk
{\it Case 2}: $\bar s=\infty$.
  Take an $s$ such that
$$
   \bar p\equiv \pi_0 \pi_1\dots \pi_s > \underline p,
$$
where $\underline p = \underline p(1/2)$ comes from Theorem \ref{improving mod-2}.
By this theorem, for sufficiently big $i$,
$$
   \mu_{n(i)- s - 1}(i) < \frac{1}{2} \mu_{n(i)}(i),
$$
contradicting again assumption (P2).
\end{proof}

We are ready to prove the Main Theorem, in an important refined version.
We say that a family $\MM$ of little Mandelbrot copies (and the corresponding
renormalization combinatorics) has {\it beau}%
\footnote{According to Dennis Sullivan,
   ``beau'' stands for ``bounded and eventually universal''. }
{\it a priori} bounds
if there exists an $\eps =\eps(\MM)>0$ and a function $N: \R_+\ra \N$
with the following property.
Let $f: U\ra V$ be a quadratic-like map with $\mod(V\sm U)\geq \de>0$
that is  at least $N=N(\de)$ times renormalizable. Then for any $n\geq N$,
the $n$-fold renormalization of $f$ can be represented  by a quadratic-like map $R^n f: U_n\ra V_n$ with
$\mod(V_n\sm U_n)\geq \eps$.

\begin{bbrmt}
  For any parameters $(\bar\q, \bar\n)$,
the family of renormalization combinatorics of decoration class $(\bar\q, \bar\n)$ has
beau a priori bounds.
\end{bbrmt}

\begin{proof}
   Let $\eps>0$ and $l$ come from Corollary \ref{improving-3},
and $C>0$ comes from  Corollary \ref{slow change}.
We will use notation $\mu_n(f)$ from the proof of Corollary \ref{improving-3}.
Assume that for some $\de>0$, there is a sequence of $\psi$-ql maps $f_i$ in question with $\mu_0(f_i)\geq \de$,
while $\mu_{n(i)}(f_i)<\eps$, where $n(i)\to \infty$. Let $n(i)=k_i l+r_i$ where $0\leq r_i<l$.
Then by  Corollaries \ref{improving-3} and \ref{slow change},
$$
   \mu_0(f_i)\leq C^l \eps /2^{r_i} \to 0\quad {\mathrm as}\ i\to \infty.
$$
This contradiction proves the beau bounds for the moduli $\mu_n(f)$ of $\psi$-ql maps.
The beau bounds for ordinary quadratic-like maps follow by Lemma \ref{ql extension}.
\end{proof}

\section{Appendix: Extremal length and width}


Given a family curves $\GG$ on a Riemann surface $S$ and a conformal metric $\mu$ on $S$,
we let $\mu(\gamma)$ be the $\mu$-length of a curve $\gamma\in \Gamma$,
$\mu(\Gamma)$ be the infimum of these lengths,
$\area_\mu$ be the corresponding measure, and
 $\LL(\GG)$ and $\WW(\GG)=\LL(\GG)^{-1}$ be respectively the {\it extremal length} and {\it width} of $\GG$:
see \cite{A} or the Appendices \cite{covering lemma,K} for the precise definitions.
The most basic properties of these conformal invariants, the {\it Parallel} and {\it Series Laws} can
also be found in these sources.

\subsection{Transformation rules}
 Both extremal length and extremal width  are conformal invariants.
 More generally, we have:

\begin{lem}\label{increase}
   Let $f: U\ra V$ be a holomorphic map between two Riemann surfaces,
and let $\GG$ be a family of curves on $U$. Then
$$
      \LL(f(\GG))\geq \LL(\GG).
$$
\end{lem}

 See Lemma 4.1 of \cite{covering lemma} for a proof.

\begin{cor}\label{two families}
     Under the circumstances of the previous lemma,
let $\HH$ be a family of curves in $V$ satisfying the following lifting
property: any curve $\gamma\in \HH$ contains an arc that lifts to some curve
in $\GG$. Then $\LL(\HH)\geq \LL(\GG)$.
\end{cor}

 See Corollary 10.3 of \cite{K} for a proof.


\msk
 Given a compact subset $K\subset \inter U$, the {\it extremal distance}
               $$\LL(U, K)\equiv \mod (U,K)$$
(between $\di U$ and $K$) is defined as $\LL(\GG)$, where $\GG$ is the family of curves
connecting $\di U$ and $K$. In case when $U$ is a topological disk and $K$ is connected,
we obtain the usual modulus $\mod(U\sm K)$ of the annulus $U\sm K$.
%
%
 We let $\WW(U,K) = \LL^{-1}(U,K)$.




\begin{lem}\label{modulus transform-2}
 Let $f: U\ra V$ be a branched covering between two compact Riemann surfaces with boundary.
Let $A$ be an archipelago in $U$,  $B=f(A)$, and
assume that  $f:  A\ra B$ is a branched covering of degree $d$. Then
$$
      \mod(V, B) \geq d\, \mod (U,  A).
$$
\end{lem}

See Lemma 4.3 of \cite{covering lemma} for a proof.


\ssk



\begin{lem}\label{modulus transform-1}
Let $(U,A)$ and $(V,B)$ be as above,
and let  $f: U\sm A\ra V\sm B$ be a branched covering of degree $N$. Then
$$
      \mod(V, B)= N\, \mod (U,  A).
$$
\end{lem}

See \cite{A} for a proof.

\subsection{Strips and quadrilaterals}

\begin{lem}\label{Pi}
 Let us consider a horizontal strip $\Pi(h)$
and an interval $I=(x, x+a)\subset \R$. We view $\Pi$ as a
quadrilateral with horizontal sides $I$ and $\R+ih$.
Then
$$
      \frac{h}{2a}\leq  \mod \Pi \leq \frac{h}{a},
$$
provided $h/a\leq 1/2$ or $\mod\Pi \leq 1/4$ (for the left-hand side inequality).
\end{lem}

\begin{proof}
  By definition, $\mod\Pi$ is the extremal length of the family of curves connecting
$I$ to  $\R+ih$. This family contains the family $\GG'$ of vertical curves
in the Euclidean rectangle with horizontal sides $I$ and $I+ih$.
Hence $\LL(\GG)\leq \LL(\GG')=h/a$.

To prove the left-hand side inequality, let us consider a Euclidean rectangle
$Q$ with vertices $x-h, x+a+h, x+a+h+ih, x-h+ih$ endowed with the Euclidean metric $\mu$.
Any curve of $\GG$ has $\mu$-length at least $h$. Hence
$$
  \LL(\GG)\geq \frac{h^2}{\area_\mu(Q)}= \frac{t}{1+2t}, \quad where\ t=h/a.
$$
We see that $\LL(\GG)\geq t/2$ for $t\leq 1/2$,
 while $\LL(\GG) >  1/4$ otherwise. The conclusion follows
\end{proof}

\begin{lem}\label{C}
Let $\Pi$ and $I$ be as in the previous lemma.
Let $C= \Pi/\, l\,\Z $
be a cylinder covered by $\Pi$
so that $I$ is embedded into the bottom of $C$.
Then
$$
      \mod \Pi  \geq \frac{1}{2} \min(\mod C,\, 0.5).
$$
\end{lem}

\begin{proof}
Since the covering $\Pi\ra C$ is an embedding on $I$, we have: $a\leq l$.
 Then by the previous lemma we obtain:
$$
    \mod \Pi \geq \frac{h}{2a} \geq \frac {h}{2l}= \frac{1}{2} \mod C.
$$
\end{proof}

\subsection{Holomorphic and embedded annuli}

Let $S$ be a hyperbolic Riemann surface with boundary with a preferred component $\si$ of $\di S$.
We assume that $S$ has finite topological type and is not the punctured disk.
A {\it holomorphic annulus}  in $\AAA$ is a holomorphic map $A: \A(1,r) \ra S$
that extends to a homeomorphism  $\phi: \T\ra \si$.
We let $\mod A = \mod \A(1,r)$.

The family of holomorphic annuli contains a subfamily of {\it embedded annuli}.
Among embedded annuli,  there is an  annulus $A_*$ of maximal modulus,
which has nice special properties. Namely, let us uniformize $A_*$ by a flat cylinder
$C= \Pi(h)/ \Z$. Then the quadratic differential $dz^2$ on $C$ is the pull-back of some
quadratic differential $q$ on $S$. Moreover, the uniformization $C\ra A_*$ extends continuously
to the upper boundary $C^+= \R+ih /\Z$ of $C$ (minus finitely many points corresponding to the punctures
of $S$),
and induces there an equivalence relation
$\tau_k : \alpha_k \ra \alpha_k'$, where $(\tau_k)$ is a finite family of  isometries between pairs of
disjoint arcs in $C^+$. The images of these arcs, $\la_k = i(\alpha_k)=i(\alpha_k')$,
are {\it horizontal separatricies} of $q$.
(It is a version of Strebel's Theorem, see e.g., \cite[\S 11]{GL}).

\begin{lem}\label{emb vs hol}
  For any holomorphic annulus $A: \A(1,r)\ra S$, we have:
$$
   \mod A\leq 16 \mod A_*.
$$
\end{lem}

\begin{proof}
  Let us consider a family $\GG$ of non-trivial proper curves $\gamma$ in $S$ that begin in  $\si$.%
\footnote{``Non-trivial'' means that $\gamma$ cannot be pulled to $\si$ through a continuous
family of proper curves.} Then  any curve $\gamma\in \GG$ contains an initial segment $\gamma'$
that lifts to a vertical curve in $\A(1,r)$. By Corollary \ref{two families},
\begin{equation}\label{upper bound}
  \mod A\leq \LL(\GG).
\end{equation}

Let us now take any conformal metric $\mu$ on $S$, and let $l= \mu(\GG)$.
For any vertical curve $\de$ in $A_*$, two possibilities can occur:

\ssk\nin $\bullet$ $\de$ ends on $\di S$. Then $\de \in \GG$ and hence $\mu(\de)\geq l$.

\ssk\nin $\bullet$  $\de$ ends on some separatrix $i(\alpha_k)$.
Then there is another vertical curve $\la$
in $A_*$ that ends at the same point as $\de$.
The concatenation of $\de$ and $\la$ is a curve of family $\GG$.
Hence one of the curves, $\de$ or $\la$, is ``long'',
i.e., it has $\mu$-length at least $l/2$.
\ssk

It follows that at least one half of the vertical curves in $A_*$ are long.
Let $I\subset C^+$ be the set of endpoints of $i^{-1}$(long curves).
We can now proceed as in the classical Gr\"otztsch estimate.
By the Cauchy-Schwarz Inequality,
$$
 h\, \area_\mu(C) = \area(C) \int_C \area_\mu\, dx\, dy   \geq  \left( \int_I dx \int_{0}^h \mu(x,y) dy\right)^2
  \geq \left( \frac{l}{4} \right)^2,
$$
which implies
$$
  \LL_\mu(\GG)=l^2/\area_\mu(C) \leq 16 h = 16 \mod A_*.
$$
Since this is valid for any conformal metric $\mu$, we conclude that $\LL(\GG)\leq 16\mod A_*$.
Together with (\ref{upper bound}), this gives us the desired estimate.
\end{proof}

\end{document}